\documentclass[review]{elsarticle}

\usepackage{float} 
\usepackage{hyperref}
\usepackage{anyfontsize}

\journal{Journal of \LaTeX\ Templates}

\usepackage{amsmath}
\usepackage{amssymb}

\newcommand{\R}{\mathbb{R}}
\newcommand{\T}{\mathbb{T}}
\newcommand{\F}{\mathcal{F}}

\begin{document}

\begin{frontmatter}

\title{Regions without invariant tori of given class for the planar circular restricted three-body problem}
\author[warwick]{N.Kallinikos}
\cortext[mycorrespondingauthor]{Corresponding author}
\ead{nikos.kallinikos@warwick.ac.uk}
\author[warwick]{R.S.MacKay}
\address[warwick]{Mathematics Institute, University of Warwick, Coventry CV4 7AL, UK}
\author[esa]{T.Syndercombe}
\address[esa]{Telespazio Germany GmbH, based at ESA/ESOC, 64293 Darmstadt, Germany}

\begin{abstract}
A method to establish regions of phase space through which pass no invariant tori transverse to a given direction field is applied to the planar circular restricted three-body problem. Implications for the location of stable orbits for planets around a binary star are deduced. It is expected that lessons learnt from this problem will be useful for applications of the method to other contexts such as flux surfaces for magnetic fields, guiding centre motion in magnetic fields, and classical models of chemical reaction dynamics.
\end{abstract}


\end{frontmatter}


\section{Introduction}

KAM theory provides sufficient conditions for existence of invariant tori in Hamiltonian systems.  Although great advances have been made \cite{FHL}, it is still hard work to obtain a realistic fraction of the tori suggested to exist by numerical simulation.

On the other hand, Converse KAM theory \cite{MP,M89}, which provides sufficient conditions for non-existence of invariant tori of given class through given regions, is much easier to implement, and in examples treated so far it produces a close to believed optimal result without much work.  Furthermore, it is proved to obtain an arbitrarily large fraction of the complement of the union of all such invariant tori under suitable conditions \cite{S}.

Converse KAM theory was developed initially for area-preserving twist maps and restricted to non-existence of invariant circles in the form of graphs:~momentum $p$ as a function of position $q$ \cite{MP} (following earlier uses by Mather, Herman and Lazutkin).  That was extended to symplectic twist maps, restricted to non-existence of invariant tori in the form of Lagrangian  graphs:~$p=\frac{\partial S}{\partial q}$ for some function $S(q)$, in general multivalued  \cite{MMS}.  Then it was extended to continuous-time Hamiltonians on $T^*\T^d$ with possibly periodic time-dependence, positive-definite second derivative in the momenta $p$ and non-existence of invariant Lagrangian graphs, and applied in particular to a Hamiltonian for the 1D motion of a particle in the field of two waves \cite{M89}.

In \cite{M18}, however, the method was extended to systems without any twist condition and for invariant tori transverse to an arbitrary foliation, in the case of 1.5 degree of freedom (DoF) Hamiltonian systems (including the restriction of an autonomous 2DoF system to energy levels). A test of the method has recently been carried out on the two-wave Hamiltonian and a quasiperiodic Hamiltonian, with encouraging results \cite{DM}.  

In this paper we apply the method of \cite{M18} to the more challenging problem of the planar circular restricted three-body problem (PCR3BP).  This venerable system is of intrinsic, practical and pedagogical interest.  Recall that it concerns the motion of a test particle in the gravitational field produced by two bodies in circular orbits around their centre of mass.  We denote the masses of the two bodies by $1-\mu$ and $\mu$ (with $\mu \in [0,\frac12]$), respectively, relative to their total mass.  We call them the primary and secondary, respectively.  The test particle is assumed to start in the plane of rotation of the two bodies with velocity in that plane and therefore to remain in that plane.  We view the motion of the test particle in a frame which keeps the centre of mass at the origin and rotates around it with the two bodies so as to keep them at positions $(-\mu,0)$ and $(1-\mu,0)$ respectively, relative to their separation.  It conserves the Jacobi constant (to be recalled in Section 3), so is a 1-parameter family of 1.5 DoF systems.

The questions we would like to address are:
\begin{enumerate}
\item What is the set of initial conditions for which the test particle is constrained to an invariant torus lying outside the orbit of the secondary ($r>1-\mu$), circulating around the origin?
\item What is the set of initial conditions for which the test particle is constrained to an invariant torus lying inside the orbit of the secondary, circulating around the primary?
\item What is the set of initial conditions for which the test particle is constrained to an invariant torus circulating around just the secondary?  
\end{enumerate}
We propose to tackle these questions by establishing the complements of the specified sets of initial conditions.  In this paper, we address just the first question, but the method could be adapted to address the other two.  Note that one could also ask about invariant tori which cross the orbit of the secondary, but continuations of these from the unperturbed case can not exist because they would include collision with the secondary.

An answer to the first question is relevant to the question of location of stable orbits for a planet round a binary star, highly topical in this age of exoplanet discovery \cite{MQ}, given that a significant fraction of stars are actually binary (estimated at around $80\%$ according to various websites).  Our paper gives initial insights into this question.

The second question is relevant to the issue of stable orbits of a small planet like the earth in the Sun-Jupiter system.

Note that an answer to the third question would provide a rigorous notion of the ``sphere of influence'' of the secondary, which has a range of definitions that scale like different powers of $\mu$, so are not compatible with each other.  There are the Hill or Roche sphere that has the line between the Lagrange points $L_1$ and $L_2$ as diameter (so scales like $\mu^{1/3}$), several variants of a sphere of influence whose radius scales like $\mu^{2/5}$
\cite{Roy}, and Belbruno's weak stability radius that scales like $\mu^{1/3}$ again \cite{CBB}.  It would be good to address this.  The three questions really require extension to the problem of motion in 3D and to the case of elliptical motion of the primary and secondary, which a planned extension of the method to higher DoF (restricting attention to Lagrangian invariant tori transverse to a Lagrangian foliation) will be able to tackle.

The plan of the paper is that first we give a simple illustration of the method.  Next we recall the Hamiltonian formulation of the PCR3BP and its invariant tori for $\mu=0$.  We state the non-existence criterion and explain how to use it in this problem.  We propose reduction of the search space of initial conditions to a surface of section or even a symmetry plane.  We give examples of initial conditions for which the method yields non-existence and then a scan of two symmetry planes to summarise results found up to a specified time-out.  We interpret the results in terms of crossing the orbit of the secondary and of resonance with the rotation of the two bodies.  We close with a discussion of various improvements that it would be good to make.

\section{Simple illustration of the method}
\begin{figure}[H]
\centering
\includegraphics[width=3.2in]{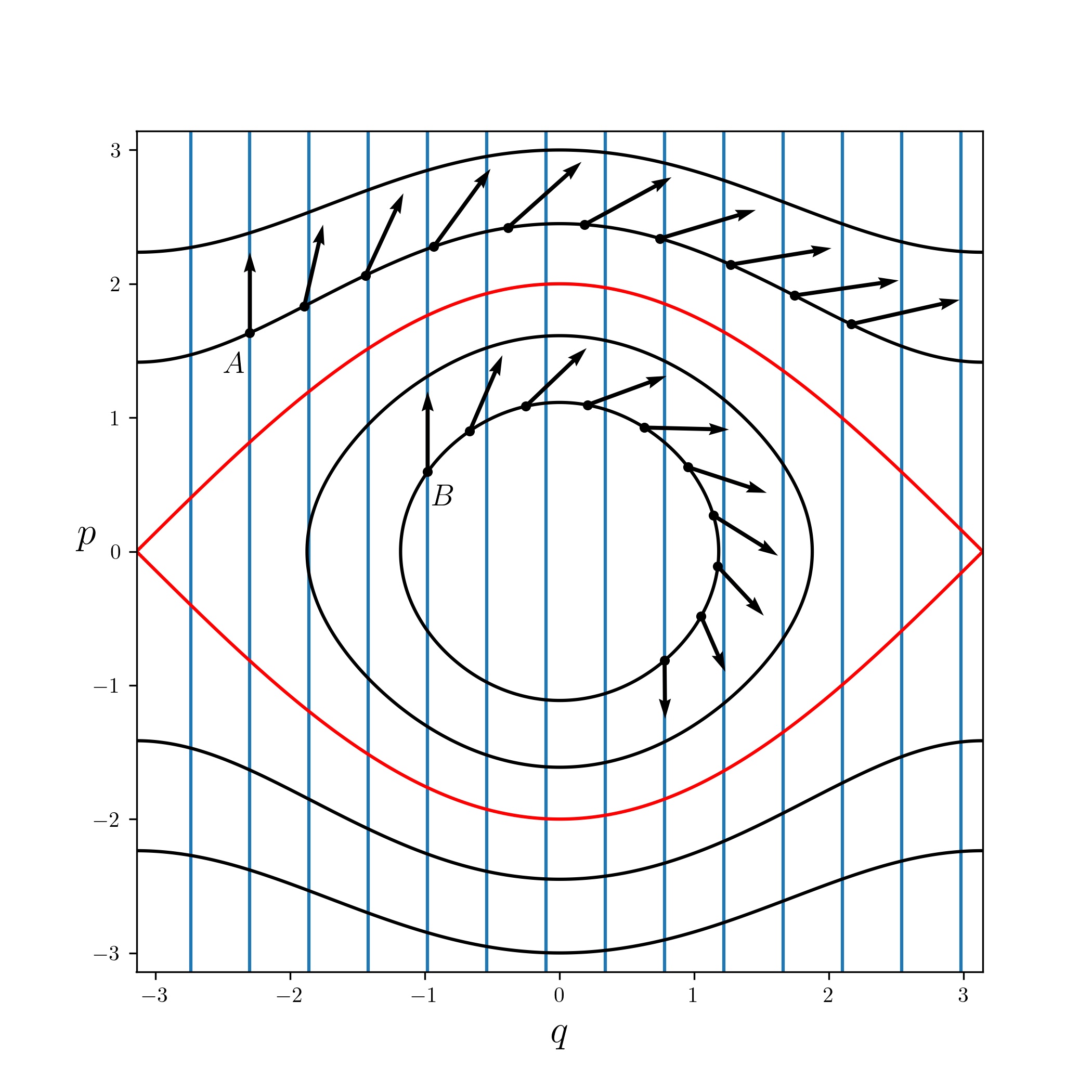}
\caption{Two vectors under the flow of the pendulum, starting at $A$ and $B$ tangent to the vertical foliation $\F$ (blue lines).}
\label{fig:pendulum}
\end{figure}

Consider the simple pendulum described by the vector field $V=(\dot{q},\dot{p})=(p,-\sin q)$ with Hamiltonian $H(q,p)=\frac12 p^2 - \cos q$ on the cylinder.  $H$ is conserved.  For $H>1$ the solutions lie on rotational invariant tori (``rotational'' means they encircle the cylinder, and here the tori are 1D, so just circles, in fact periodic orbits).  They are all transverse\footnote{Two submanifolds are said to be \textit{transverse} if at any intersection the sum of their tangent spaces is the whole tangent space.} to the foliation\footnote{A \textit{foliation} of a manifold is a decomposition into subsets called \textit{leaves} that are locally submanifolds, diffeomorphic to the decomposition of $\R^n$ into the set of leaves $x=$ constant for some choice of $0<m<n$ and coordinates $(x,y) \in \R^{n-m} \times \R^{m}$.}
$\F$ given by the vertical lines $q=$ constant. Given a foliation $\F$, we define a vector field $\xi$ by a continuous choice of upward tangents to $\F$. The orbit $\eta$ of an upward tangent $\xi$ to $\F$ at a point $A$ with $H > 1$ under the linearised dynamics cannot cross the tangent to the invariant torus, therefore it cannot become a downward tangent to $\F$, as illustrated in Figure~\ref{fig:pendulum}. Thus, if an upward (or downward) tangent to $\F$ at some point flows to a downward (or upward) tangent, then no rotational invariant tori pass through the given point (nor any point of its orbit).

For $H<1$ therefore, we can tell that the solutions are not on rotational invariant tori, by taking the orbit of an upward tangent to $\F$ at any such point $B$ and noticing that it crosses the downward tangent after some time, as in Figure~\ref{fig:pendulum}.

Here, of course, the solutions for $H<1$ lie on librational invariant tori. In order to test nonexistence for this type of torus, we could choose instead a radial foliation (e.g.~$p/q=$ constant, though this does not extend to a global foliation of the cylinder), and one would find for this example that the test is never satisfied.

One could rightly say that an even simpler condition, and one that would be satisfied earlier, is that if a trajectory lies on a torus of the given class then $\eta$ never becomes parallel to the dynamical vector field $V$.  But the above formulation allows extension to higher dimensions.  In particular, we extend the above condition slightly to encompass both of the above ideas. If a trajectory is on an invariant torus of the given class then the trajectory $\eta$ of $\xi$ from an initial point can never become of the form a negative amount of $\xi$ plus an arbitrary amount of $V$, because that would imply that $\eta$ becomes tangent to the torus at some point in between, but the trajectory of a tangent to the torus remains forever tangent to it.

\section{PCR3BP}
We choose units in which the total mass is 1, the distance between the primary and secondary is 1, and the angular frequency of rotation is 1.  We choose coordinates $(x,y)$ for the test particle in the rotating frame so that the primary mass $1-\mu$ is at $(-\mu,0)$ and the secondary mass is at $(1-\mu,0)$.

Let ${\bf p}=(p_x,p_y), L, K$ be the vector momentum, angular momentum and energy (per unit mass) of the test particle in the instantaneous inertial frame, and ${\bf r} = (x,y)$ be its position.  So
\begin{gather}
L = xp_y - y p_x\\
K = \frac12\left(p_x^2+p_y^2\right) - \frac{1-\mu}{r_1}-\frac{\mu}{r_2},
\end{gather}
where $r_1$ and $r_2$ are the distances to the primary and secondary, respectively, $r_1=\sqrt{(x+\mu)^2+y^2}$ and $r_2 = \sqrt{(x-1+\mu)^2+y^2}$.

Then the motion in the rotating frame is given by Hamilton's equations with respect to the canonical symplectic form $\omega =  dx \wedge dp_x + dy \wedge dp_y$ for 
\begin{equation}
H = K-L.
\end{equation}
We denote the resulting vector field by $V$ (satisfying $i_V \omega = dH$). In particular, $H$ is conserved and it is conventional to denote its value by $-C/2$, with $C$ called the Jacobi constant.

\begin{figure}[H]
\centering
\includegraphics[width=3.8in]{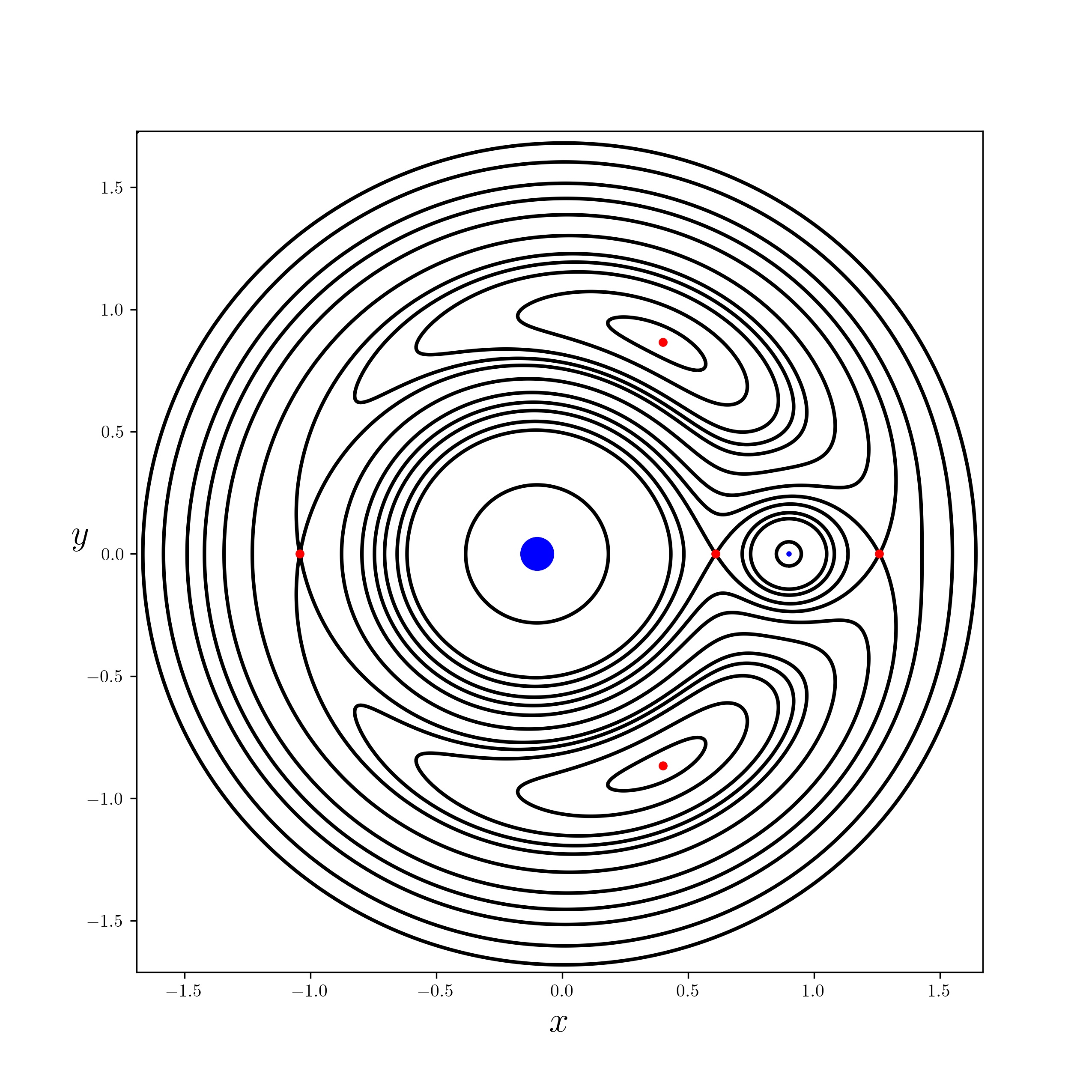}
\caption{Contours of the effective potential $U$ for $\mu=0.1$. The two massive bodies are shown in blue, and the five Lagrange points are shown in red. $U$ goes to $-\infty$ at the two masses and has maxima at the equilateral Lagrange points.}
\label{fig:potential}
\end{figure}

The system can alternatively be written in terms of the velocity ${\bf v}=(v_x,v_y)$ in the rotating frame.  The transformation is
\begin{align}
\begin{split}
v_x &= p_x+y,\\
v_y &= p_y-x.
\end{split}
\end{align}
The Hamiltonian and the symplectic form become
\begin{gather}
H = \frac12\left(v_x^2+v_y^2\right)+U(x,y),\\
\omega = dx \wedge dv_x + dy \wedge dv_y - 2 dx \wedge dy
\end{gather}
with
\begin{equation}
U(x,y) = -\frac12\left(x^2+y^2\right) - \frac{1-\mu}{r_1}-\frac{\mu}{r_2}.
\end{equation}
Figure \ref{fig:potential} shows a contour plot of the effective potential $U$.  Its critical points express the equilibrium points of the system, its level sets $U=-C/2$ are called ``zero-velocity curves'' and its sub-level sets $U\le -C/2$ are the ``Hill's regions'' for allowed motion at Jacobi constant $C$.

For state space, we use $S = \{(x,y,v_x,v_y) \in \R^4: (x,y) \ne (-\mu,0), (1-\mu,0)\}$ and endow it with Euclidean metric (using ${\bf v}$ rather than ${\bf p}$), which we will use to construct the foliation vector field $\xi$. It is possible to regularise the system to allow passage through collisions in a modified time, but we leave that for the future.

\subsection{The unperturbed system}

\begin{figure}[H]
\centering
\includegraphics[width=3in]{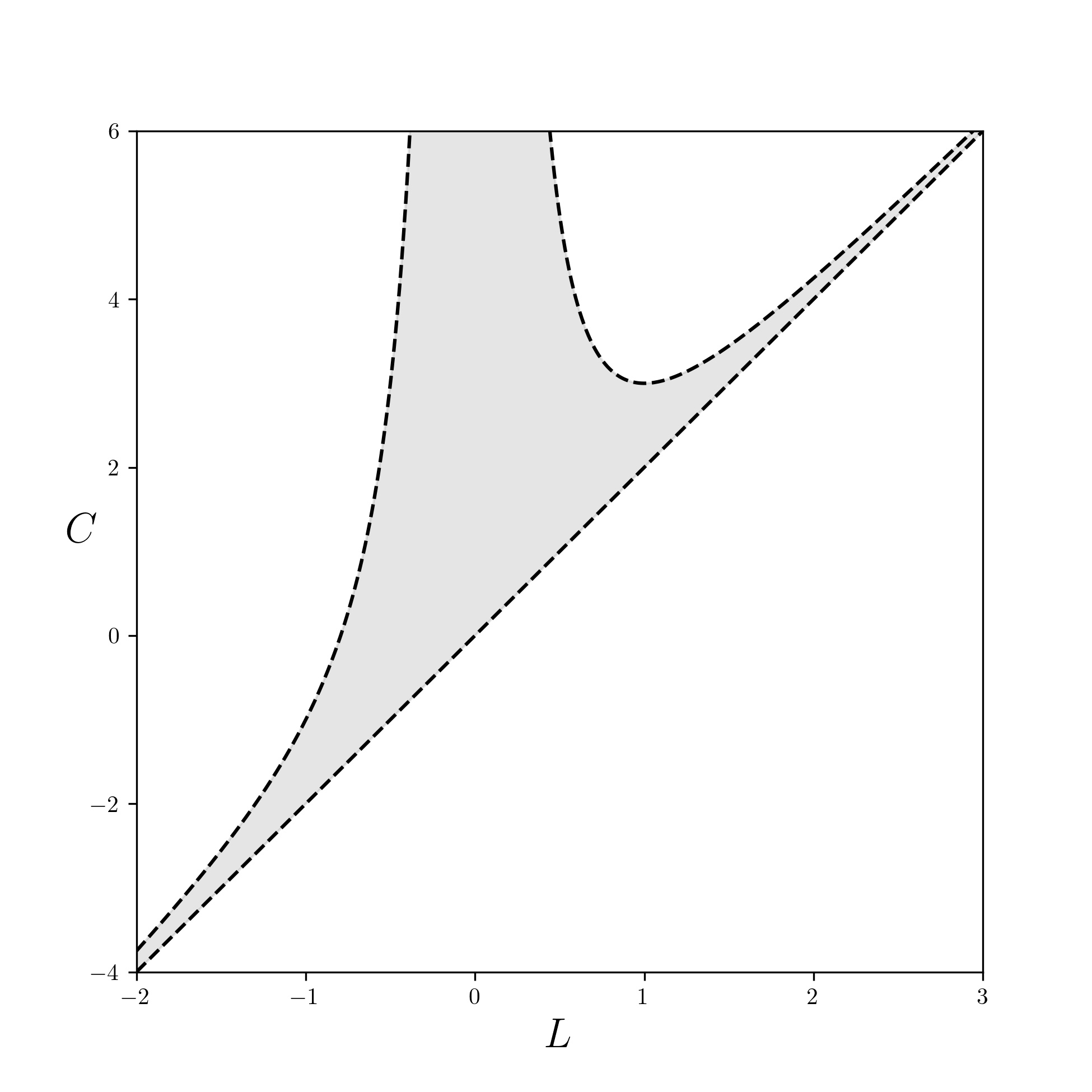}
\caption{The region of invariant tori (grey) in the space of $(L,C)$ for $\mu=0$.}
\label{fig:CLregion}
\end{figure}

When $\mu=0$, the system is integrable, with first integrals $K$ and $L$.  The joint level sets of $(K,L)$ are invariant and connected.  They are non-empty, bounded and regular  iff $K<0$ and $0< L^2 < (-2K)^{-1}$ (``regular'' means that the derivatives of $K$ and $L$ are linearly independent everywhere on them).  The non-empty bounded regular level sets of $(K,L)$ are two-tori. In polar canonical coordinates $(r,\theta,p_r,p_{\theta})$, where $L=p_{\theta}$, the joint level sets of $(K,L)$ can be equally described by
\begin{equation}
\label{eq:invariants}
p_r^2+\frac{L^2}{r^2}-2L-\frac{2}{r} = -C.
\end{equation}
The two-tori correspond to parameters $2L<C<2L+L^{-2}$, shown in Figure~\ref{fig:CLregion}. The region is bounded above by curves corresponding to the circular orbits $p_r=0, r=L^2$, and below by a line corresponding to the parabolic orbits.

Each torus corresponds to a choice of semi-major axis $a>0$, eccentricity $e\in(0,1)$ and direction $\sigma\in\{\pm 1\}$ of rotation (direct or retrograde).  The relation is $K=-(2a)^{-1}$, $L=\sigma \sqrt{a(1-e^2)}$.  The points on a torus with given $a,e,\sigma$ correspond to points on the Kepler ellipses with these parameters; the remaining freedoms are the angle $g$ of pericentre in $(x,y)$ (closest approach to the origin) and the position of the particle on the ellipse.  The latter can be described by the polar angle $\theta$ or the angle $f$ relative to pericentre or the mean anomaly $m$ (being $2\pi$ times the fraction of the area swept out from pericentre) or by the eccentric anomaly $E$ (that we won't use) \cite{Roy}. The dynamics on such a torus are conjugate to a constant vector field:
\begin{equation}
\dot{m}=N^{-3}, \quad \dot{g}=-1,
\end{equation}
where $N = \sigma\sqrt{a}$ is the first Delaunay variable (so $K=-(2N^2)^{-1}$).  The equation $\dot{g}=-1$ comes from the rotating frame. The winding ratio of turns in $g$ to turns in $m$ is $w = -N^3$.

\begin{figure}[H]
\centering
\includegraphics[width=2.39in]{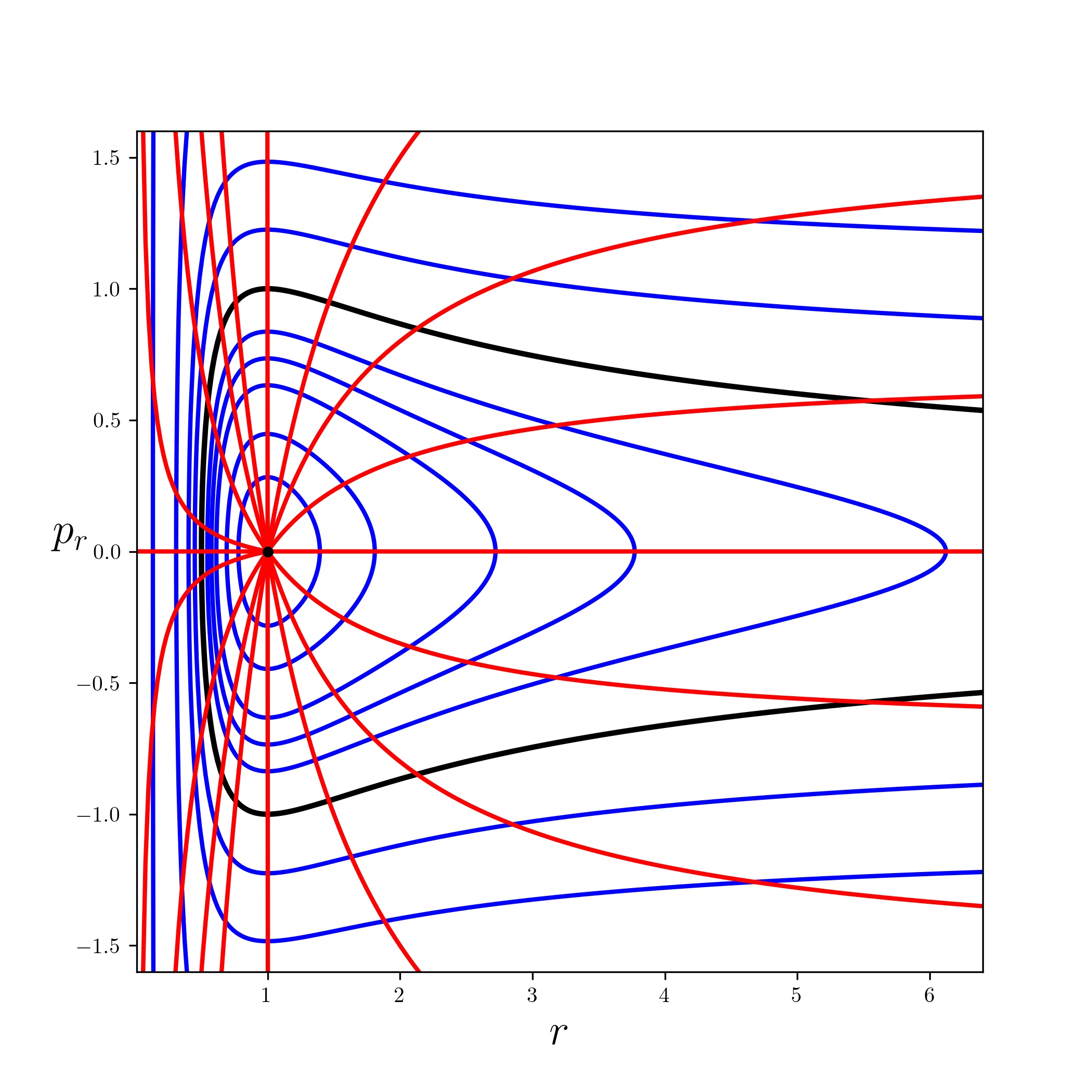}\includegraphics[width=2.39in]{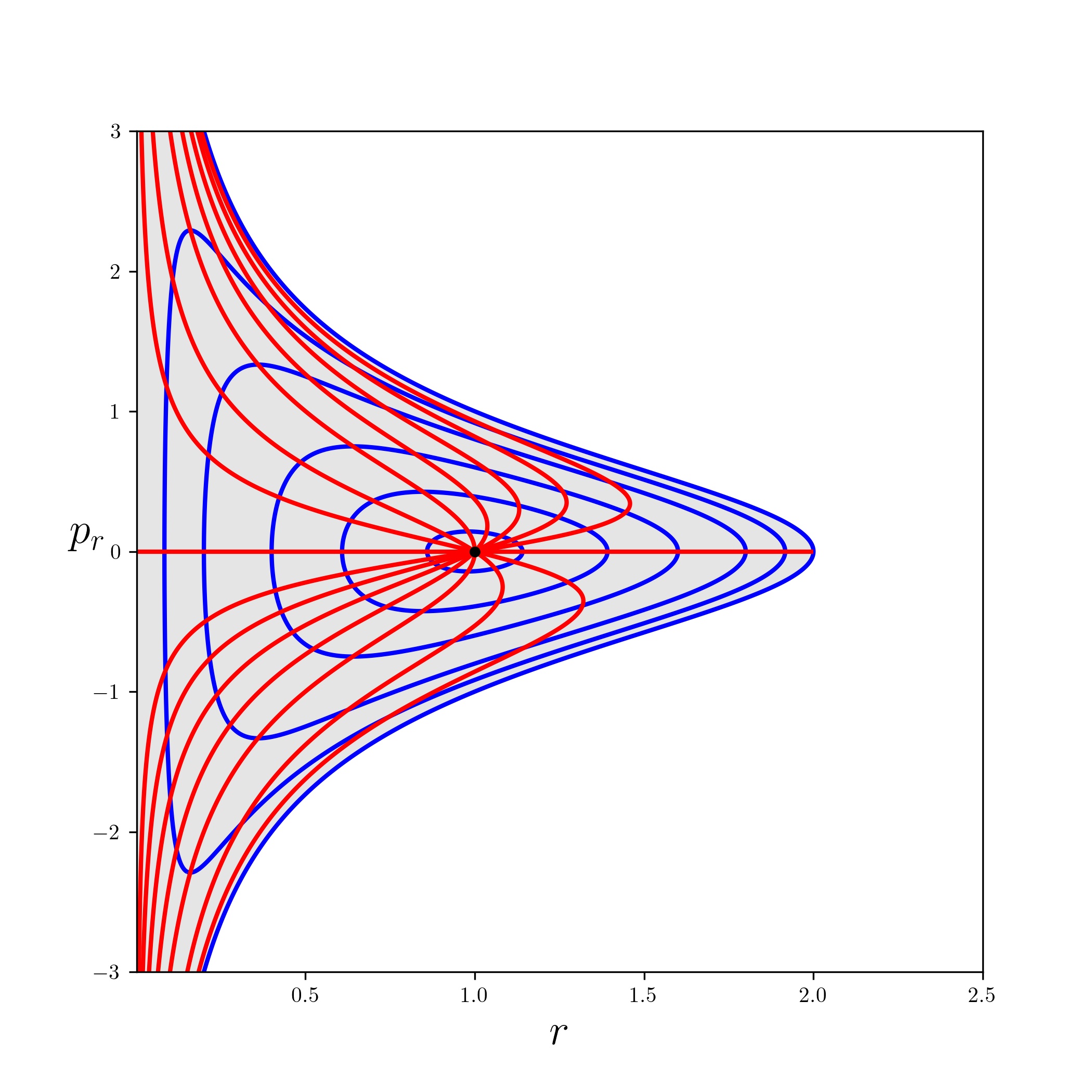}
\caption{Joint level sets of $(K,L)$ (blue) and transverse foliation $\F$ (red) for $\theta=$ constant, in the unperturbed case. For fixed $L=1$ (left), the tori range from $K=-0.5$ for circular orbits (black point) to $K=0$ for parabolic ones (black curve). For fixed $K=-0.5$ (right), the allowed motion (grey area) is confined by the outermost tori $L=0$ and degenerates to the circular orbits for $L=\pm1$.}
\label{fig:invariants}
\end{figure}

The invariant surfaces are transverse to the foliation $\F$ defined by $g=$ constant and $\theta=$ constant.  This is because $(L,N,g,\theta)$ forms a local coordinate system (except at $r=0$ where $\theta$ is undefined, and on the circular orbits where $g$ is undefined). The foliation $\F$ has singularities corresponding to the coordinate singularities, so we should call it a singular foliation. Back to polar coordinates, the 2D leaves of $\F$ are given by the level sets of $F=p_rL(L^2/r-1)^{-1}$. This is because $g$ is the angle of the Laplace vector (often called Runge-Lenz) ${\bf e} = {\bf p} \times L {\bf\hat{z}}- {\bf\hat{r}}$ to the positive $x$-axis. The intersections of some of the invariant surfaces and leaves of the foliation $\F$ with $\theta=$ constant and either $L=1$ or $K=-0.5$ are shown in Figure~\ref{fig:invariants}.\vspace{-0.05cm}

\section{Converse KAM method}

By KAM theory, sufficiently irrational invariant tori persist smoothly and thus remain transverse to $\F$ for some range $0\le \mu < \mu_c(w,L,K)$, where $\mu_c$ depends smoothly on $L,K$, but in a fractal number-theoretic manner on $w$ (for a heuristic description of this in the context of the two-wave Hamiltonian, in terms of what are now called Brjuno functions, see \cite{M88}).  For question 1, it may be more appropriate to consider $\mu$ as fixed and express the conclusion of KAM theory as applying when the distance $r_{\min}= a(1-e) = N^2(1-\sqrt{1-L^2/N^2})$ to pericentre of the osculating ellipse satisfies $r_{\min}  > r_c$ for some function $r_c(\mu,w,L)$ depending smoothly on $\mu$ and $L$, but number-theoretically on $w$.

The Converse KAM method of \cite{M18} eliminates regions for a 3D system where invariant tori transverse to a given foliation don't exist. So, to apply the method we first need to choose a foliation for the energy levels of the PCR3BP.

\subsection{Foliation}

A natural candidate is based on the foliation $\mathcal{F}$ introduced at the end of the previous section. Restricting to $H=-C/2$, $\F$ becomes 1D.  This turns out not to be a good choice of foliation of the energy level for $\mu>0$, however, because the leaves develop tangencies to the energy levels when the perturbation $\mu$ is turned on, introducing extra singularities and making it difficult to introduce a consistent orientation, a feature needed by the non-existence method.  The issue is that the effect of the perturbation is not small near the primary and secondary. To answer question 1, it might seem that we do not have to consider trajectories that come close to the primary or secondary, but in practice it turned out that we needed a way to handle them, so we decided it is tidier to choose a foliation of the energy levels that does not have this problem.

To specify a 1D oriented foliation for each energy level, it is enough to specify a vector field $\xi$ tangent to the energy levels.  Based on the idea that the invariant tori for $\mu=0$ in an energy level are $L=$ constant, we chose 
\begin{equation}
\label{eq:xi}
\xi = \nabla L - a \nabla H, 
\end{equation}
with $a = \nabla L \cdot\nabla H /|\nabla H |^2$, where $\nabla, \cdot, |\,|$ are with respect to the standard Riemannian metric on $(x,y,v_x,v_y) \in \R^4$. The metric  mixes lengths and velocities, which might seem physically unsatisfactory, but we have already scaled lengths and times to make the distance between the primary and secondary and the rotation rate be one. The vector field $\xi$ is undefined where $\nabla H = 0$, but that is only the five Lagrange points, where in any case the energy levels have singularities. By construction, $\xi$ is tangent to the energy levels ($\xi \cdot \nabla H = 0$), and it is transverse to the level sets of $L$ ($\xi \cdot \nabla L >0$) except where $\nabla L$ and $\nabla H$ are parallel. In the unperturbed case $\mu=0$, these are parallel only at $r=0$ and on the circular orbits $L^2=r$, $p_r = 0$.  Thus $\xi$ is transverse to the invariant 2-tori of the unperturbed case ($L=$ constant in the domain of Figure~\ref{fig:CLregion}) in the given energy level.

The vector field $\xi$ induces an oriented foliation of each energy level by its integral curves.  The foliation has singularities where $\xi=0$, i.e.~where $\nabla L$ and $\nabla H$ are parallel.  For $\mu=0$ this was already discussed.  For $\mu>0$ it consists of two curves in the full phase space.  They are most conveniently written using $\bf p$ rather than $\bf v$, and they can be expressed as $p_r=0$, $\dot{p}_r=0$, $\dot{p}_{\theta}=0$ in polar coordinates.  More specifically, the first one is given by $\theta=0$, $p_r=0$, $f(r,L;\mu)=0$ and $\theta=\pi$, $p_r=0$, $f(r,L;1-\mu)=0$, where
\begin{equation}
\label{eq:sigbound}
f(r,L;\mu)=\frac{1-\mu}{(r+\mu)^2} \pm \frac{\mu}{(r-1+\mu)^2} - \frac{L^2}{r^3},\quad \pm(r-1+\mu)>0.
\end{equation}
This is a deformation of the equation for the circular orbits, restricted to $y=0$; the deformation is small except near $x = 1-\mu$ or $-\mu$. It will be illustrated in Figure~\ref{fig:symlevels}.
The second one is $r_1=r_2$, $p_r=0$, $L^2 = r^4 r_1^{-3}$,
which again is a deformation of the equation for circular orbits, but restricted to the perpendicular bisector of the massive bodies.

\subsection{Nonexistence condition}\label{subsec:nonexistence}

Now we present from \cite{M18} a sufficient condition for non-existence of invariant tori through a given point, transverse to the vector field $\xi$, adapted here to the PCR3BP.  

Take an initial point $s_0$ in $H^{-1}(-C/2)$ and an initial tangent vector $\eta_{s_0} = \xi_{s_0}$. For increasing $t$, simultaneously evolve both to $s=s(t)$ and $\eta_s=\eta_s(t)$ using the dynamics $\dot{s}=V(s)$ and the linearised dynamics $\dot{\eta}_s = DV_s\,\eta_s$. If there is an invariant torus $\mathcal{T}$ passing through $s_0$ that is transverse to $\xi$, then $\eta$, that is, $\eta_{s(t)}$ for all $t$, must stay on the same side of $\mathcal{T}$. In particular, we can never have (i) $\eta_s,\xi_s,V_s$ linearly dependent, with (ii) $\eta_s = \alpha V_s + \beta \xi_s$, $\beta < 0$.

To detect (i), we can use the symplectic form $\omega$ for the system, because in a regular energy level, $\eta_s,\xi_s,V_s$ are linearly independent if and only if $\omega(\eta_s,\xi_s)\ne 0$. To see this, take the triple product on $H^{-1}(-C/2)$ using the Liouville volume-form $\Omega = |V|^{-2}V^\flat\wedge \omega$.  Here, $V^\flat$ is the 1-form such that $V^\flat(X) = V\cdot X$ for all vectors $X$, and $\Omega$ is the standard volume-form on an energy surface inherited from Liouville volume on the whole state space $S$ (such that $\Omega \wedge dH = \frac12 \omega\wedge\omega$). So we look for a sign change of $\omega(\eta_s,\xi_s)$.

To decide (ii), we reformulate \cite{M18} by choosing a 1-form $\lambda$ such that $\lambda(V) = 0$ and $\lambda(\xi)>0$ (except at zeroes of $\xi$). So then any vector $\eta_s$ tangent to an energy level that satisfies $\lambda(\eta_s)<0$ and is dependent on $(V_s, \xi_s)$, has $\beta<0$ (in the above notation).  We say informally that $\eta_s$ points opposite to $\xi$, relative to the Hamiltonian vector field $V$.

Thus, putting the two together we arrive at the converse KAM condition: If there is a point $s=s(t)$ where $\omega(\eta_s,\xi_s)$ changes sign and $\lambda(\eta_s)<0$ then there is no invariant torus through $s_0$ transverse to $\xi$. We will refer to this as the \textit{general formulation}.

The only thing that remains is to choose $\lambda$. For the choice (\ref{eq:xi}) for $\xi$, we take $\lambda=dL-bV^\flat$ with $b=\nabla L \cdot V /|V|^2$. By construction, this satisfies $\lambda(V) = 0$. Now, note that $|\xi|^2=\xi\cdot\nabla L$ and $V\cdot\xi=V\cdot\nabla L$, because $\xi$ and $V$ are each perpendicular to $\nabla H$. Thus, $\lambda$ and $\xi$ also satisfy $\lambda(\xi)=\nabla L\cdot\xi-b\,V\cdot\xi=|\xi|^2-(V\cdot\xi)^2/|V|^2\ge0$ by the Cauchy-Schwarz inequality. Therefore $\lambda(\xi)$ is positive everywhere except where $\xi$ is parallel to $V$.

For convenience, we allow the possibility that $\xi$ is parallel to $V$ in some places. Nonexistence of invariant tori transverse to $\xi$ through such a point is automatic. The condition $\lambda(\eta_s)<0$, however, is not satisfied there as $\lambda(\eta_s)=0$ at these points. This might be unfortunate, but we did not come up with a choice of $\xi$ that we could guarantee to be nowhere parallel to $V$.

We close this section with an alternative formulation of the converse KAM method, which although not used here, might be helpful in the future. One way to choose $\lambda$ is to choose a vector field $u$ independent from $\xi$ (except at its zeroes) and tangent to the energy levels, such that $\alpha(u)>0$, where $\alpha=i_{\xi}\omega$, and then let $\lambda = -i_u \omega$. It follows that $\lambda(V)= dH(u) = 0$ and $\lambda(\xi)=\alpha(u)>0$ are automatic. Moreover, we don't need to work with $\lambda$ directly:~instead of flowing a tangent vector $\eta_s$, one can take an initial cotangent vector $\beta_{s_0}=\alpha_{s_0}$ and let it flow to $\beta_s=\beta_s(t)$ under the adjoint linearised system $\dot{\beta}_s = -\,\beta_sDV_s$. Since $\beta_s=i_{\eta_s}\omega$, we arrive at the following condition. If there is a point $s=s(t)$ where $\beta_s(\xi_s)$ changes sign and $\beta_s(u)<0$ then there is no invariant torus through $s_0$ transverse to $\xi$.

\subsection{Symmetric case}\label{subsec:symmetric}

A refinement of the nonexistence condition, which goes back to \cite{M89}, involves systems that admit a time-reversal symmetry $(t\longrightarrow-\,t,\,s\longrightarrow\tilde{s}=R(s))$ for some diffeomorphism $R$. The PCR3BP has indeed the time-reversal symmetry $R:(r,\theta,p_r,p_{\theta})\longrightarrow(r,-\theta,-p_r,p_{\theta})$. 

For time-reversal symmetric systems and initial conditions on a symmetry plane (the set of fixed points of $R$), we can speed up the computations by a factor of at least two by using the fact that the backward trajectory is the reflection by $R$ of the forward one.  Thus we get a segment of trajectory of twice the length for the price of one.  One should choose the vector field $\xi$ to be symmetric with respect to $R$, i.e. $\tilde{\xi}_{\tilde{s}} = \xi_s$ where $\tilde{\xi} = dR\, \xi$.  Instead of starting with $\eta_{s_0}=\xi_{s_0}$, we  choose any antisymmetric $\eta_{s_0}$ on the symmetry plane independent from $V_{s_0}$ (which is automatically antisymmetric), but tangent to the energy level sets.

Then the non-existence condition can be refined to $\omega(\eta_s,\xi_s)=0$ for some $t>0$.  This is because, if $\eta$ is antisymmetric, it follows that if one starts from $\tilde{s}$ at $-\,t$ with tangent vector $-\,\tilde{\eta}_{\tilde{s}}$ then one obtains tangent vector $\eta_{s_0}$ at $t=0$ and hence $\eta_s$ in $s$ at time $t$.  But $\omega(\eta_s,\xi_s)=0$ with $dH(\eta_s)=0$ implies $\eta_s= c_1 \xi_s + c_2 V_s$ for some $c_1,c_2 \in \R$, with $c_1 \ne 0$ (because the only way to get $c_1=0$ is to start with $\eta_{s_0}$ a multiple of $V_{s_0}$, but we took it independent).  Then we deduce that starting from $\tilde{\eta}_{\tilde{s}}= -\,c_1\tilde{\xi}_{\tilde{s}} + c_2\tilde{V}_{\tilde{s}}$ at $\tilde{s}$ produces $\eta_s = c_1 \xi_s+c_2 V_s$ at $s$, which is incompatible with having an invariant torus through $s_0$, transverse to $\xi$.

In conclusion, given now an $R$-symmetric $\xi$ and an $R$-antisymmetric $\eta_{s_0}$ (both tangent to the energy levels), if there is a point $s=s(t)$ where $\omega(\eta_s,\xi_s)$ changes sign then there is no invariant torus through $s_0$ transverse to $\xi$. We will refer to this as the \textit{symmetric formulation} of converse KAM.

Note that the above condition is satisfied automatically where $\xi$ becomes parallel to $V$. Therefore this formulation does not rule out invariant tori through such points, even though they cannot be transverse to $\xi$ and lie outside of the class in question, and will correctly pick up nonexistence there.

\begin{figure}[H]
\centering
\includegraphics[width=3.1in]{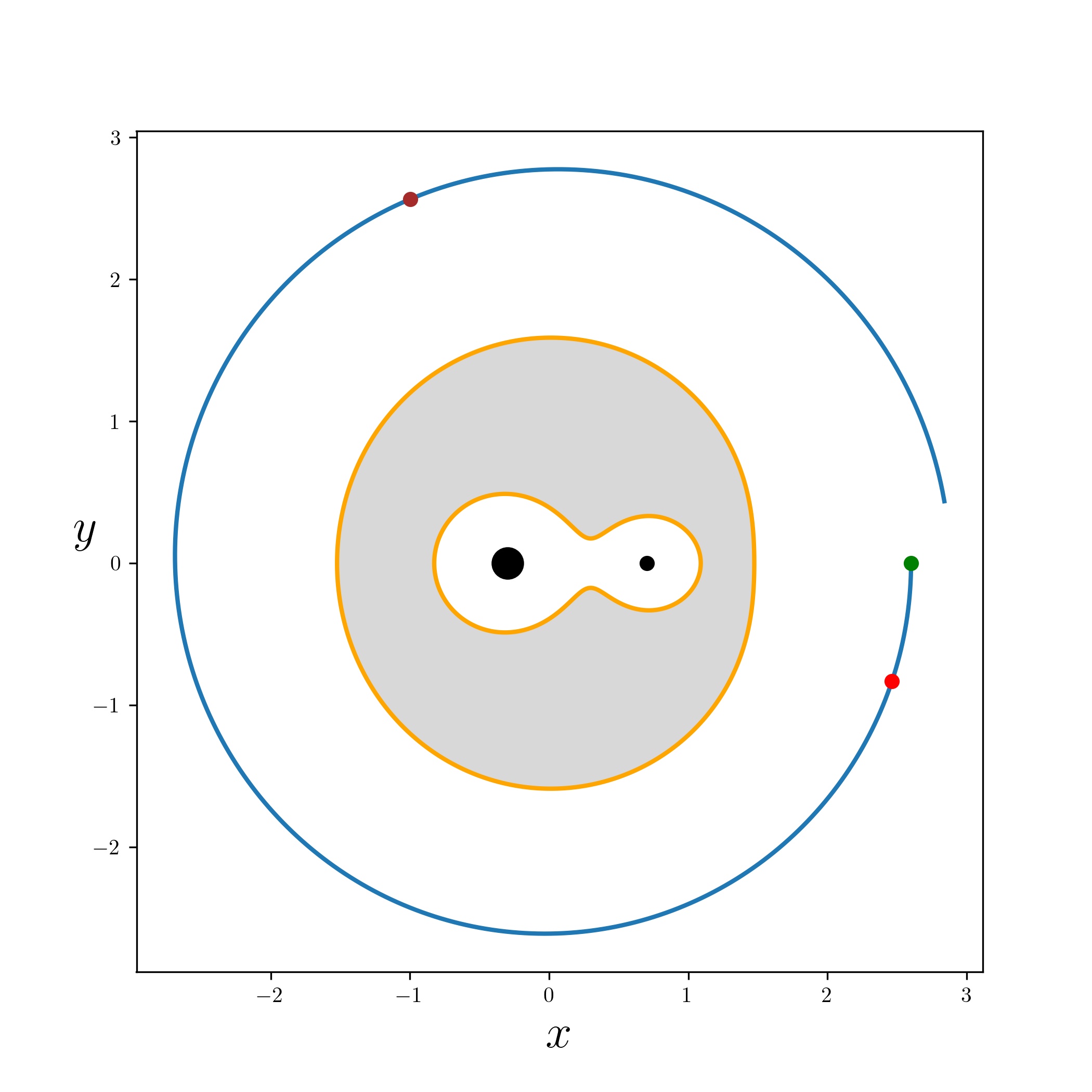}
\caption{Example of a trajectory for which the nonexistence condition was satisfied. The two bodies are shown in black, the grey area is the forbidden region bounded by the zero-velocity curves (yellow), and $\mu=0.3$, $C\approx3.7$ with initial conditions (green point) on a symmetry plane. The symmetric formulation of converse KAM detects here nonexistence (red point) much faster than the general formulation (brown point).}
\label{fig:orbit}
\end{figure}

Figure~\ref{fig:orbit} shows an example of a trajectory, which converse KAM detected in both formulations, and the general formulation of the nonexistence condition is compared to the refined one using the time-reversal symmetry.

\section{Reduction of dimension}

\subsection{Surface of section for bounded orbits}

It is enough to test initial conditions on a codimension-1 set $\Sigma$ such that every bounded trajectory crosses it.  Such a $\Sigma$ is called a surface of section if it is transverse to the vector field $V$, but this is not necessary for present purposes.

Every bounded trajectory comes to a local maximum of $r$, so take $\Sigma = \{s \in S: p_r=0, \dot{p}_r\le 0\}$ (a trajectory could have $r$ increasing to a supremum as $t \to \infty$, but that happens only for  trajectories approaching a Lagrange point).  For $\mu=0$ this is $\{s \in S: p_r=0, L^2\le r \}$.  

We can examine one value of $C$ at a time.  We denote $\Sigma_C = \Sigma \cap H^{-1}(-C/2)$, which is 2D. For $\mu=0$ the allowed region on $\Sigma$ is $2L\le C \le 2L+L^{-2}$, as was shown in Figure~\ref{fig:CLregion}, so $\Sigma_C$ consists of one or two annuli according as $C\le 3$ or $C>3$.  Unfortunately, for $\mu \ne 0$ the effect of the secondary is large near $\theta = 0$ and the effect of displacement of the primary is large near $\theta=\pi$.  The result is that there are large deviations of the allowed region from the case $\mu=0$ near these angles.

\begin{figure}[htbp]
\centering
\includegraphics[width=3.6in]{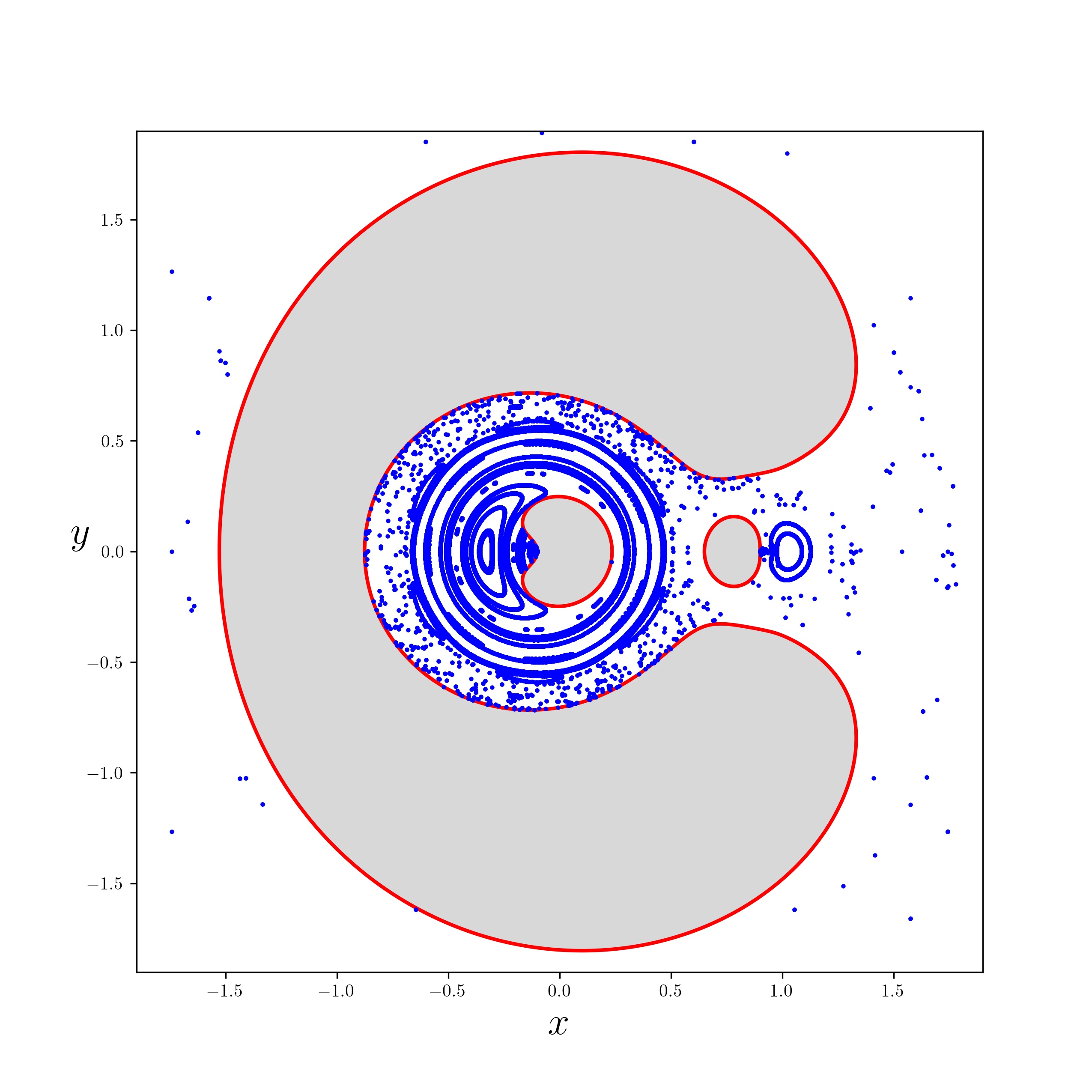}\\\vspace{-0.4cm}
\includegraphics[width=4.4in]{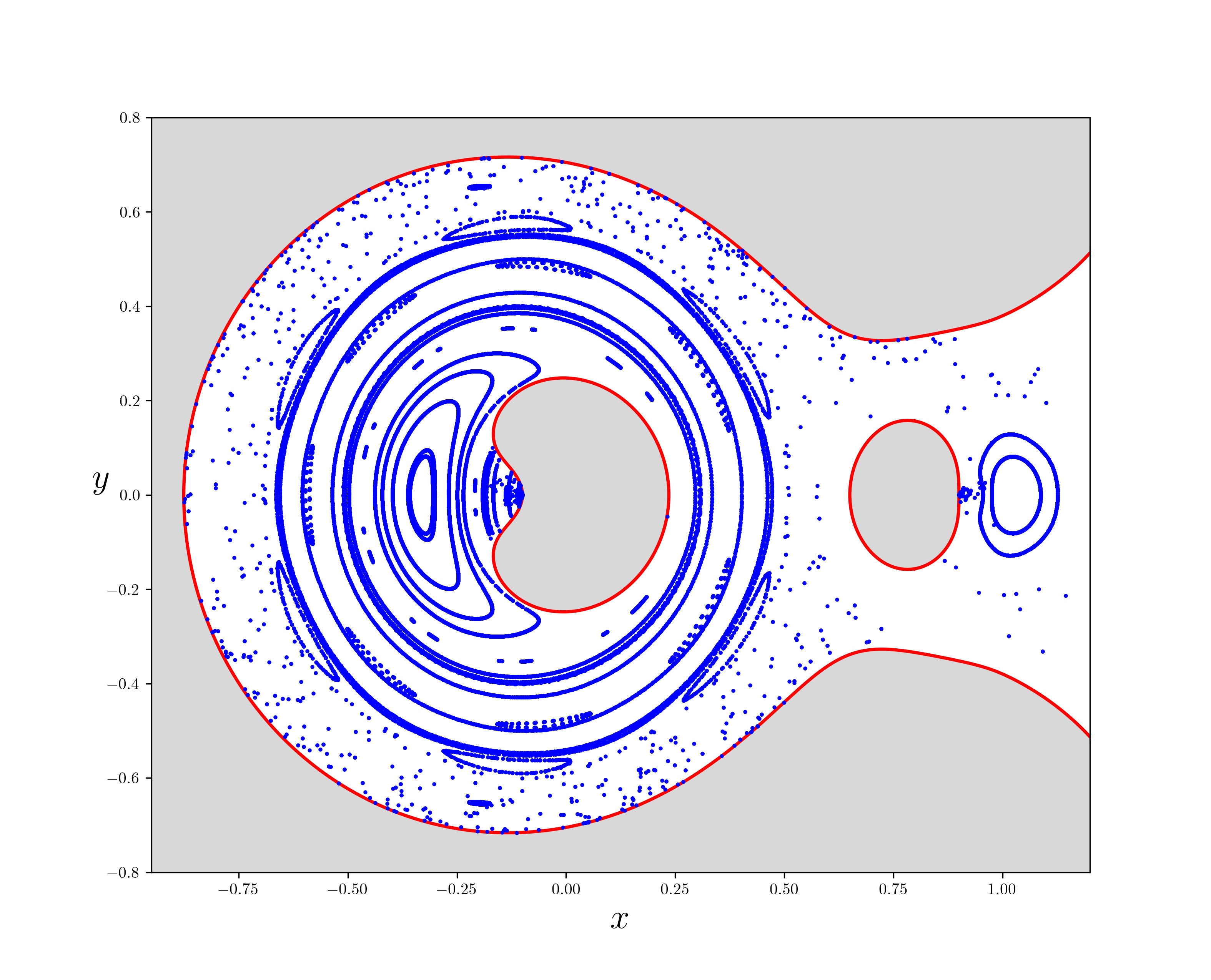}
\caption{Some orbits of the return map to $\Sigma_C$ (top), projected to the $(x,y)$-plane, for $\mu=0.1$, $C=3.2$, and zoomed in (bottom) around the two bodies. The grey area is the forbidden region for $\Sigma_C$.}
\label{fig:sos}
\end{figure}

Figure~\ref{fig:sos} shows the successive returns to $\Sigma_C$ for $\mu=0.1, C=3.2$, for some trajectories. Although for $\mu=0$, the restricted surface of section $\Sigma_C$ has natural coordinates $(r,L)$, for $\mu>0$ the 2D surface $\Sigma_C$ and therefore (the image of) the return map to it cannot be one-to-one mapped to a plane in any of the usual coordinates. Ideally, we would deform the surface of section to a Birkhoff section, that is a codimension-1 surface which is transverse to the flow except on its boundary which is invariant under the flow.  That requires finding the continuation of the circular orbits to $\mu>0$, however, so is not a straightforward prospect.  An alternative would be to adopt the approach of \cite{DW}, but we decided to do neither, because we do not really require a surface of section; it is just for illustration.

We note in passing that there is an additional sufficient condition we could use for the PCR3BP, namely that if $s(t)$ never returns to $\Sigma$ then it is not on an invariant torus (of any class).  This is because every trajectory on an invariant torus is bounded and so must come to a local maximum of $r$ (the possibility that a trajectory has $r$ increasing to a supremum as $t \to \infty$ is excluded on an invariant torus).  We already used this condition to exclude $K\ge 0$ for $\mu=0$, but it would be good to use it for $\mu>0$ because the non-existence condition of the previous section does not distinguish between invariant tori and invariant submanifolds that go to infinity.  One ought to be able to find an explicit condition that guarantees $p_r$ remains positive forever after.

We can reduce the search by one more dimension if we choose just a single leaf of the foliation in each energy level (i.e., an integral curve of $\xi$), because every 2-torus transverse to $\xi$ has to cut that leaf.  If we choose the selected leaves smoothly with respect to energy then they make a 2D surface $P$.  Thus to exclude an invariant torus transverse to $\xi$ it suffices to exclude the corresponding point on $P$.
Two catches are that we don't know which point on $P$ corresponds to a given torus and there might be points of $P$ not on invariant tori for which the non-existence condition is never satisfied.  Nevertheless, if for example, we establish that no points of $P$ are on invariant tori transverse to $\xi$ then we deduce that there are no such invariant tori.  This was used in \cite{MP} to prove that the standard map has no rotational invariant circles for any parameter value $k\ge 63/64$.

\subsection{Symmetry planes}

A particularly natural choice for $P$ is a symmetry plane with respect to a time-reversal symmetry $R$. Recall that the PCR3BP has the time-reversal symmetry $R:(r,\theta,p_r,p_{\theta})\longrightarrow(r,-\theta,-p_r,p_{\theta})$.  The symmetry planes are the sets of fixed points of $R$, namely $P_0$ and $P_\pi$ defined by $p_r=0$ and $\theta=0, \pi$ respectively.  We can use coordinates $(r,L)$ on them.  

Note that in satisfying $p_r=0$ the symmetry planes have some commonality with the surface of section.  Specifically, if we restricted attention on a symmetry plane to $\dot{p}_r\le0$ then it would be a 2D subset of the 3D surface of section.  Indeed, on a symmetry plane the part with $r\ge L^2$ for $\mu=0$ corresponds to $\dot{p}_r \le 0$. Likewise, for $\mu>0$ the part with $f\ge0$ (see (\ref{eq:sigbound})) corresponds to $\dot{p}_r \le 0$. This is because the boundary of $\Sigma$ coincides with the singularities for $\xi$ on the symmetry planes.

\begin{figure}[H]
\centering
\includegraphics[width=2.39in]{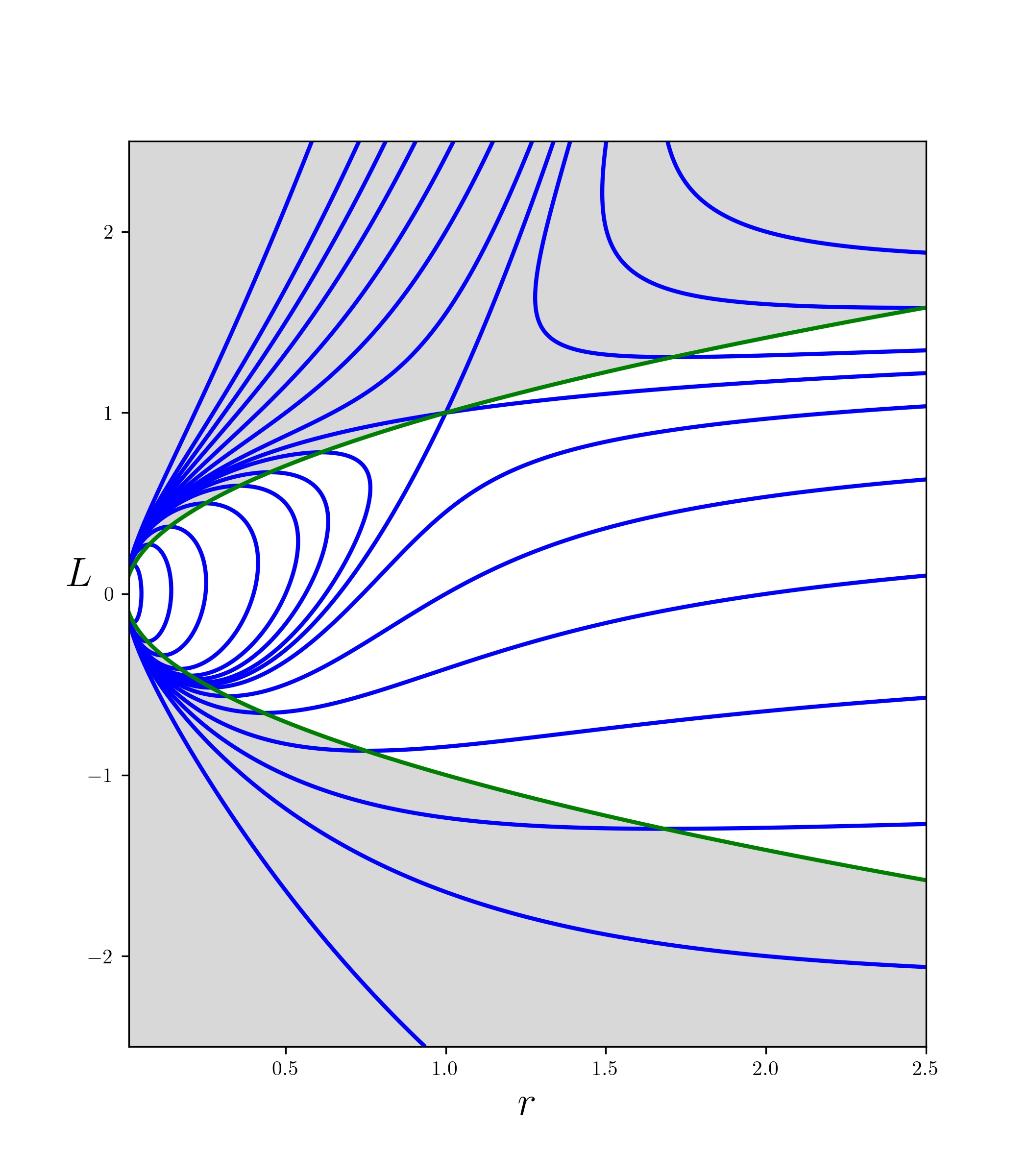}\includegraphics[width=2.39in]{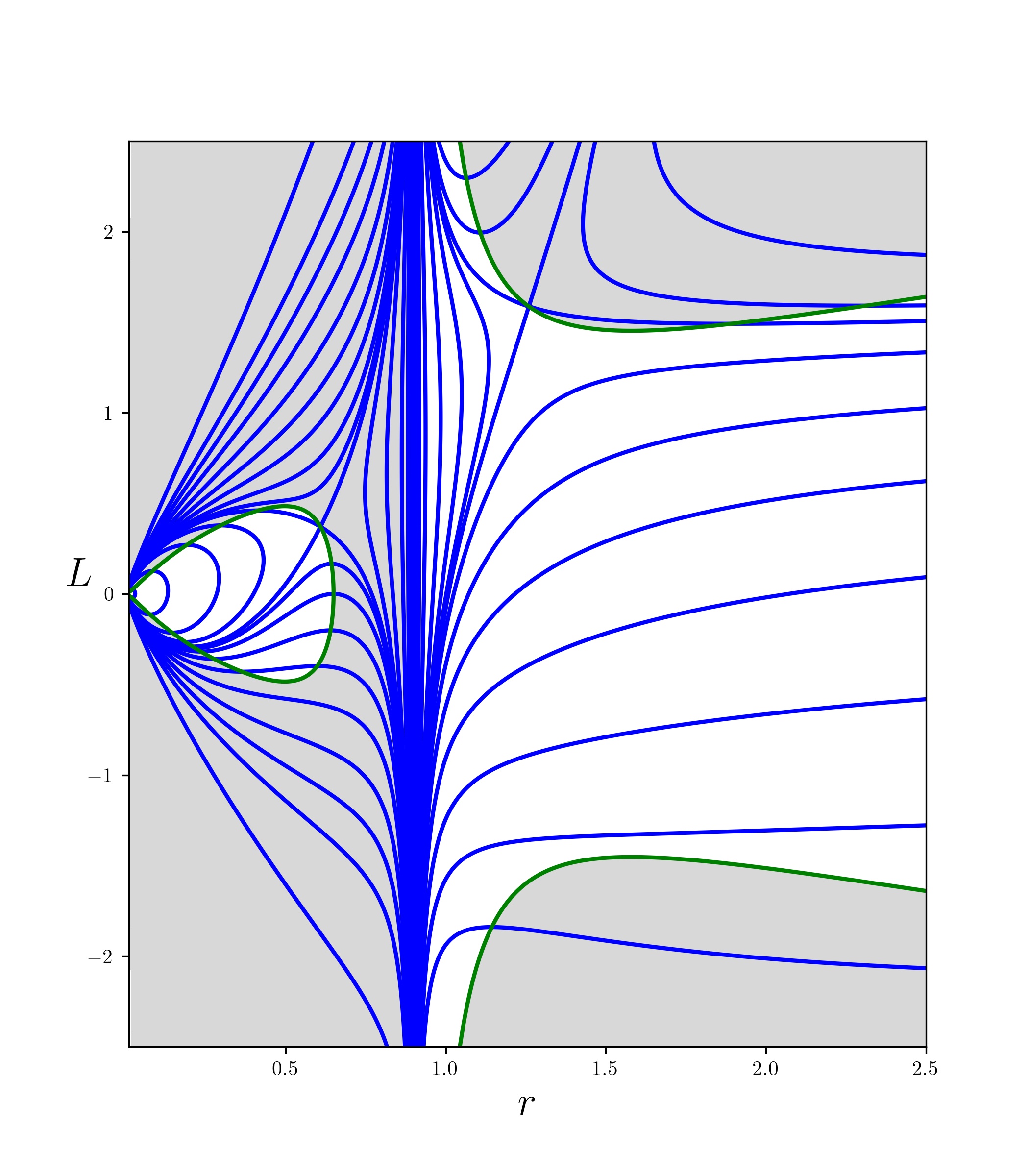}
\caption{Level sets of $H$ (blue) and the surface of section $\Sigma$ (white) bounded by $f=0$ (green) on the symmetry plane $P_0$ for $\mu=0$ (left) and $\mu=0.1$ (right).}
\label{fig:symlevels}
\end{figure}

Again, one could treat one value of $C$ at a time.  The intersection of a symmetry plane with $H^{-1}(-C/2)$ is a curve or pair of curves, as shown in Figure~\ref{fig:symlevels}.  But as the symmetry planes are only 2D, this decomposition into $C=$ constant is hardly useful. Note, one can restrict to $\dot{p}_r\le 0$ if desired too.

\section{Results}

In this section, we apply the converse KAM condition using both the general and the symmetric formulations. We present the results for initial conditions $s_0$ in the symmetry planes $P_0$, $P_{\pi}$, from which we exclude any singularities for $\xi$ when implementing the general formulation, and any singularities for $\eta$ when implementing the symmetric formulation. Following numerically both the flow and the tangent flow of the system, we track which orbits satisfy the converse KAM condition for nonexistence of invariant tori transverse to the $\xi$-direction within a fixed timeout $t_{\text{out}}$. In the plots that follow, we indicate in red initial conditions that correspond to nonexistence, and in blue that no result was obtained before timeout. The excluded singularities in each case are shown in black.

For comparison, we also include in green the resonances of the unperturbed system. The formula for resonance with rational ratio $w$ is $(w^{-2/3}r-1)^2 = 1-w^{-2/3}L^2$, giving ellipses in $(r,L)$.  Furthermore, we include in yellow the curve $L^2=2r/(r+1)$, which bounds the region ($L^2<2r/(r+1), r>1$) in which the orbit of the initial condition starts outside and crosses the orbit $r=1$ of the secondary for $\mu=0$ (to see this, use $r_{\max} = a(1+e)$, $r_{\min}=a(1-e)$ and $L^2=a(1-e^2)$ for Kepler ellipses). In grey, we also indicate the curve $L^2=2r$ which is the boundary $K=0$ of the set of initial conditions whose orbits remain bounded for $\mu=0$. Lastly, recall that the singularities for $\xi$ in black serve also as a boundary of $\Sigma$ on the symmetry planes.

\begin{figure}[H]
\centering
\includegraphics[width=4in]{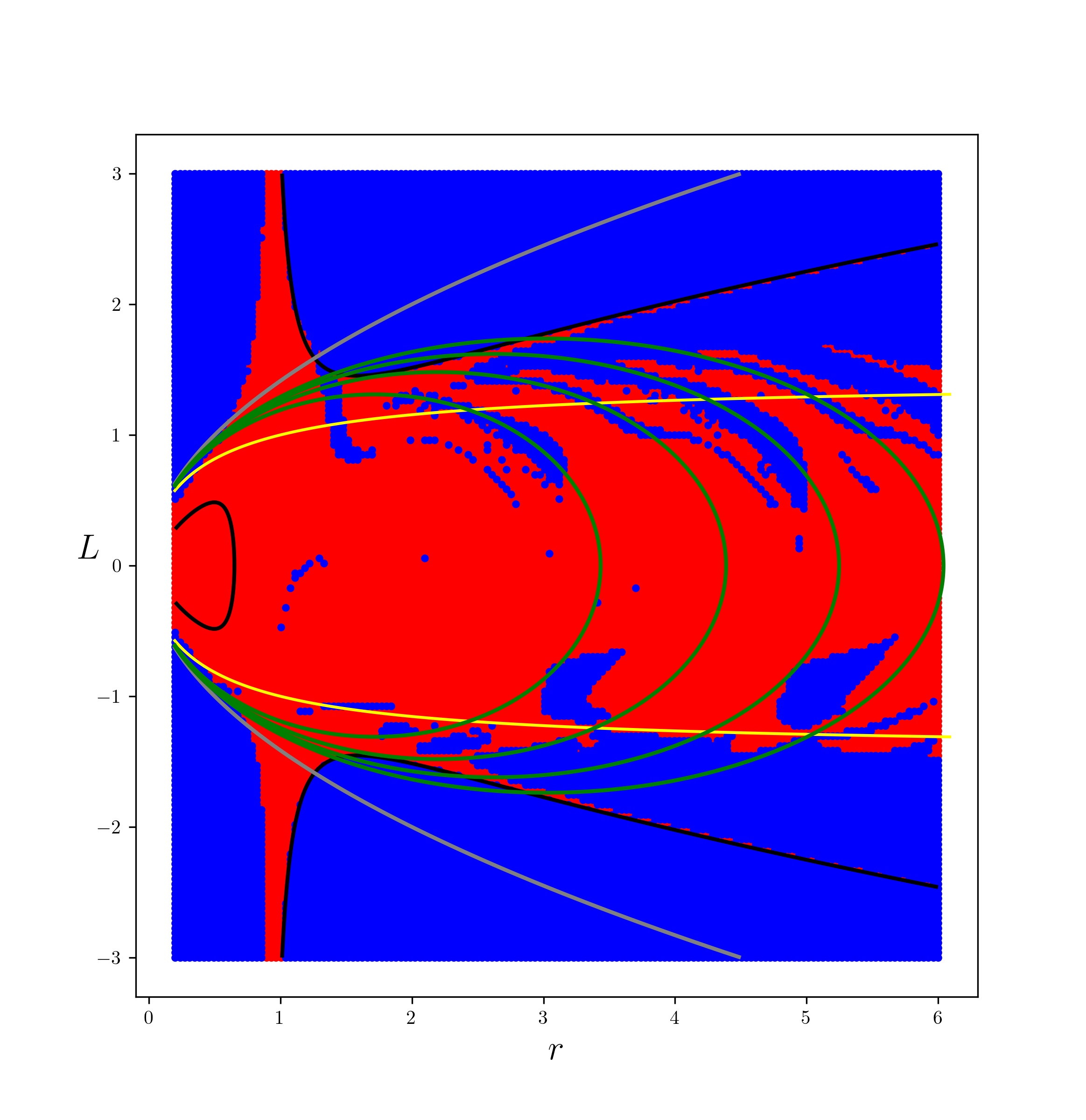}
\caption{Converse KAM results from the general formulation for $\mu=0.1$ on the symmetry plane $P_0$. Red = nonexistence, blue = undetermined, black = zeroes of $\xi$. Superposed are the resonances in green with winding ratios $n/4$ for $n=9,13,17,21$.}
\label{fig:lambda}
\end{figure}

\begin{figure}[htb]\vspace{-1.2cm}
\centering
\includegraphics[width=4in]{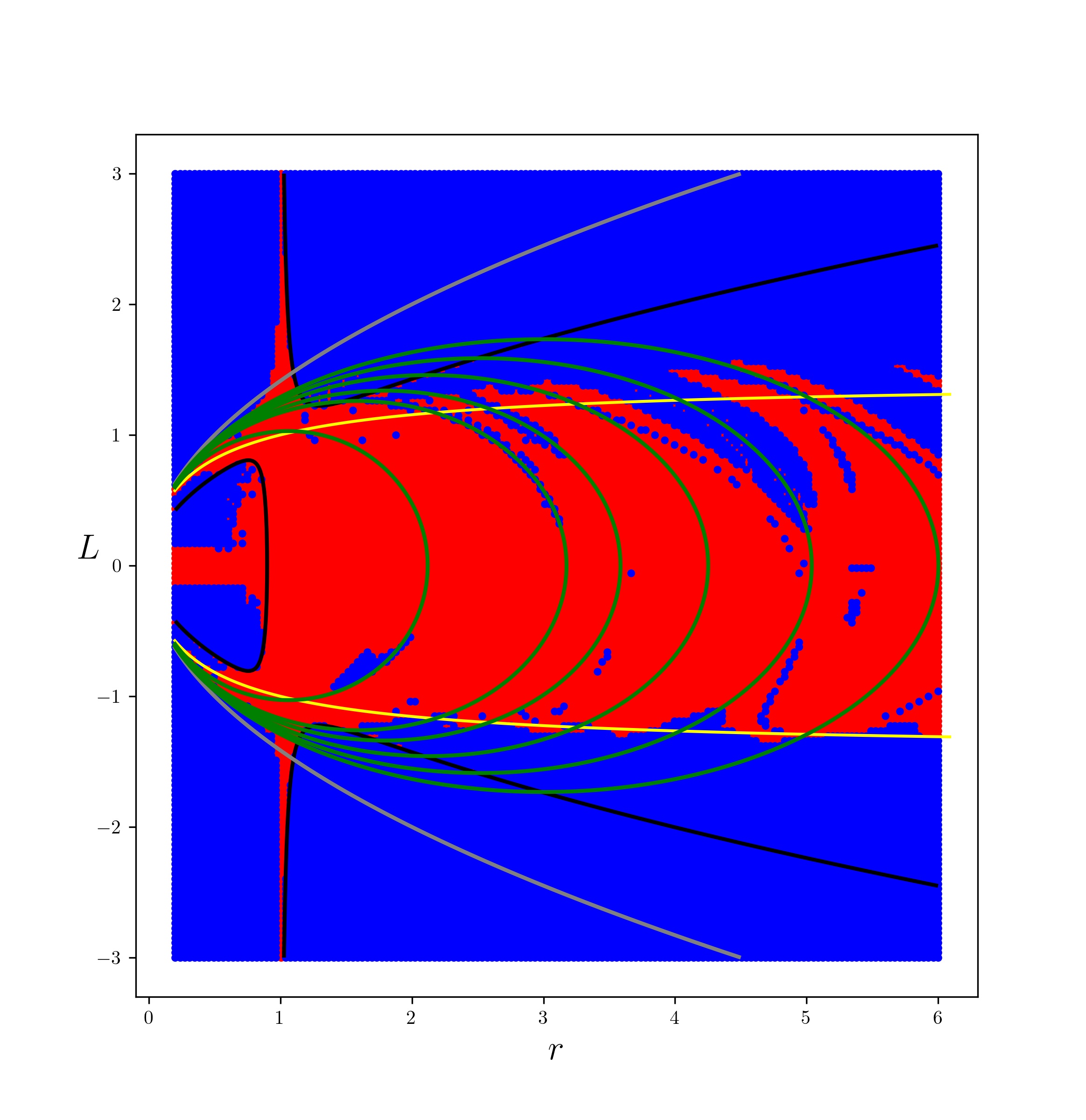}
\caption{Converse KAM results for $\mu=0.01$ (the rest of the setup same as in Figure~\ref{fig:lambda}). Superposed are the resonances in green with winding ratios $12/11,2,12/5,31/10,4,26/5$.}
\label{fig:lambda01}
\end{figure}

We start with the general formulation and the symmetry plane $P_0$, for which the results are shown in Figures~\ref{fig:lambda} and \ref{fig:lambda01} for $\mu=0.1$ and $\mu=0.01$, respectively.\vspace{-0.45cm}

\begin{figure}[htbp]
\centering
\includegraphics[width=3.6in]{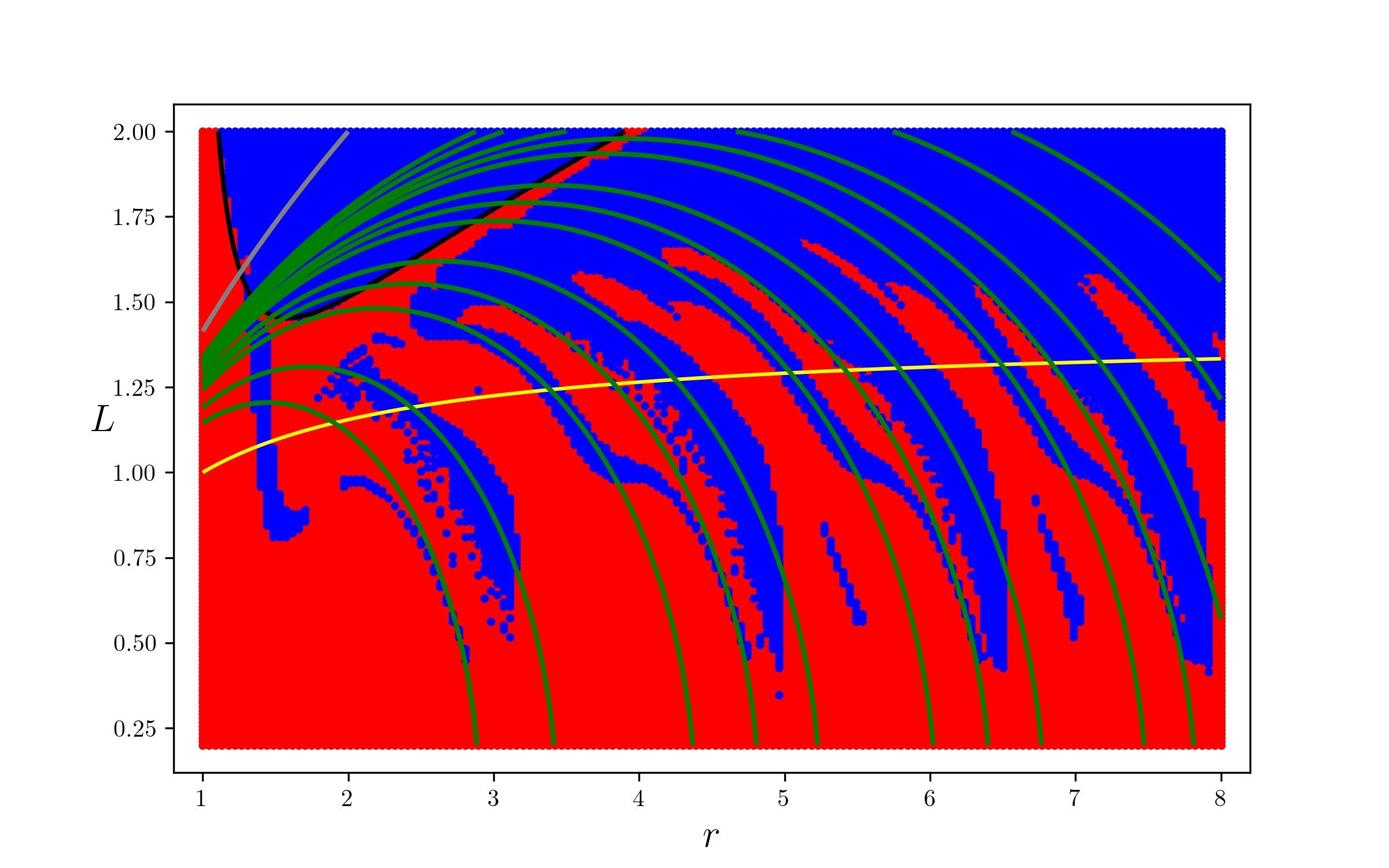}\\\vspace{-0.1cm}
\includegraphics[width=3.6in]{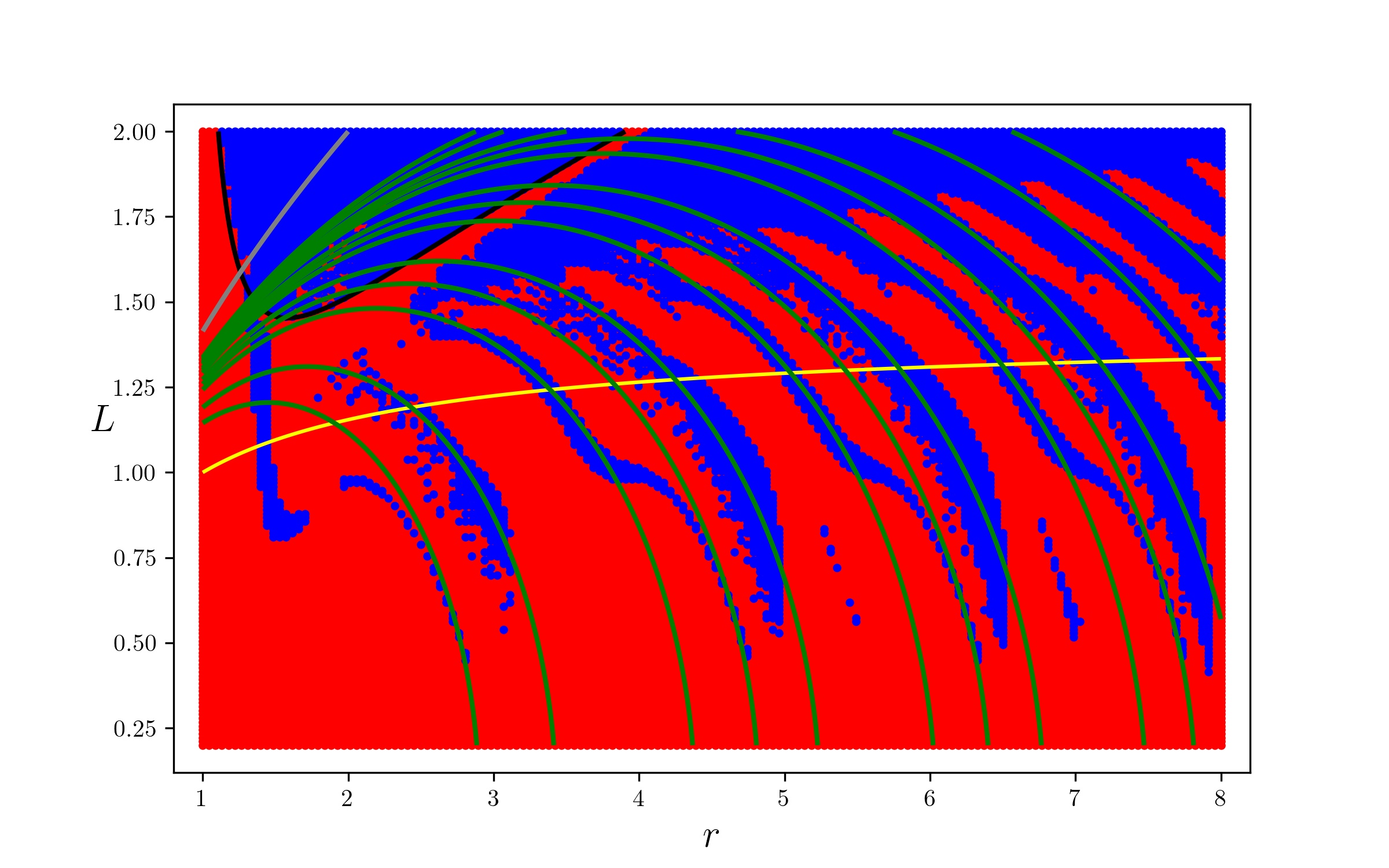}\\\vspace{-0.1cm}
\includegraphics[width=3.6in]{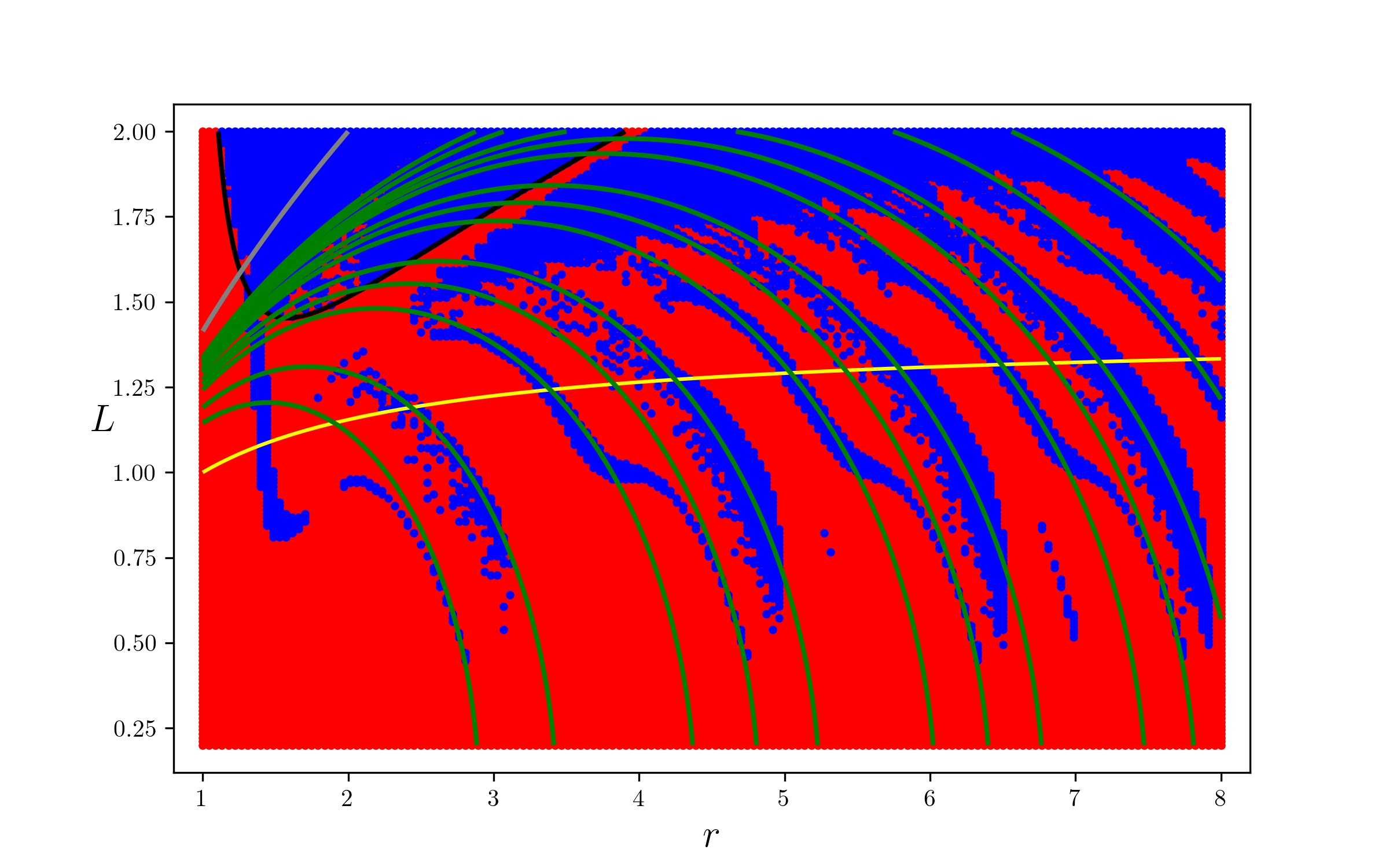}
\caption{Converse KAM results of Figure~\ref{fig:lambda} zoomed in the region of interest and run for progressively longer timeout periods (double (middle) and triple (bottom) timeout). This does appear to reach a limiting case, though longer time periods can also increase the opportunity for numerical inaccuracies to produce false positives.}
\label{fig:timeout_comparison}
\end{figure}

The code implementation is such that the numerical integration of each orbit stops if it detects the nonexistence condition or if it reaches a selected maximum timeout. Running for a longer time  results in identifying more initial conditions on which the method succeeds, and ultimately all orbits satisfying the converse KAM condition would be identified. Figure \ref{fig:timeout_comparison} shows how longer time-periods do indeed appear to demonstrate that the results tend to a limiting case.

\begin{figure}[H]
\centering
\includegraphics[width=4in]{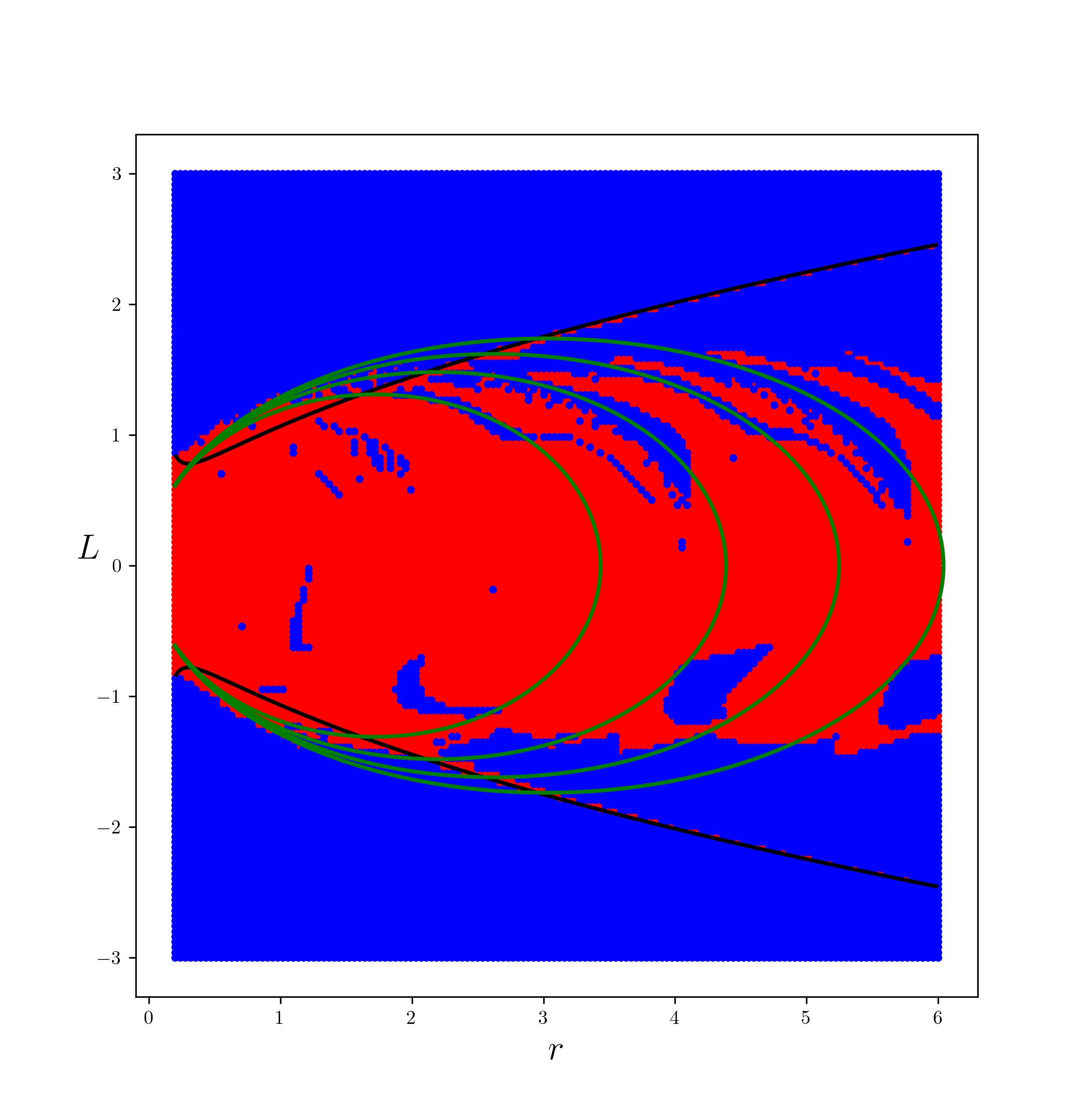} 
\caption{Converse KAM results on the symmetry plane $P_{\pi}$ (the rest of the setup same as in Figure~\ref{fig:lambda}).}
\label{fig:lambdapi}
\end{figure}

From these results and, in particular, Figure~\ref{fig:lambda}, which we will use as a reference case, we see and comment accordingly the following:
\begin{enumerate}
\item The subset with $r<1$, $f>0$ and a significant fraction of the area around it where the method gives non-existence of tori (recall that $f$ is the function defined in (\ref{eq:sigbound})). This is not surprising for orbits that come within distance $\mu$ of the primary because the foliation is based on Kepler ellipses about the centre of mass, whereas more appropriate for small or eccentric orbits would be about the primary.  Similar remarks go for orbits which approach the secondary.  It is more surprising for the rest of $r<1, f>0$, but perhaps $\mu=0.1$ is large enough that all trajectories feel so much influence from the primary and the secondary that the tori are destroyed.
\item A strip around $r=1$ where the method gives non-existence of invariant tori of the given class. Initial conditions with $r$ near $1$ in $P_0$ are close to the secondary, so it is not surprising that their dynamics lead to something very different from an invariant torus of the given class.
\item A large fraction of the subset with $r>1$ and $L^2<\frac{2r}{r+1}$ is shown to have no invariant tori of the given class. When $\mu=0$, all initial conditions in this subset produce Kepler ellipses that cross the orbit $r=1$ of the secondary. Indeed, this is precisely the condition that $r_{\text{min}}<1$ for given initial $r=r_{\text{max}}$. It is not surprising that after some time depending on their relative frequencies, the test particle should suffer a near collision with the secondary and that this should lead to a trajectory very different from an invariant torus of the given class.  Exceptions are initial conditions for which a resonance maintains a positive minimum distance from the secondary.
\item There are some points near $f=0$, $r>1$ for which the method shows no invariant tori of the given class. The set $f=0$ corresponds to points where the foliation is singular, and for $\mu=0$, in particular, to circular orbits. Although these orbits are surrounded by invariant tori and many of these are expected to persist as $\mu$ increases, they will in general deform and the thinnest ones will fail to be transverse to the foliation because of its singularity.
\item Low-order resonance for some rationals produces a significant zone where the method gives non-existence, but not for all rationals. We would expect resonance for $\mu=0$, where $w$ is a low-order rational, would lead to zones of non-existence of invariant tori of the given class for $\mu>0$ because of the formation of island chains.  We were surprised at first to see this for only some rationals.  A glance at Figure~\ref{fig:sos}, however, shows that on $\theta=0$ the principal island chains all have hyperbolic points.  A feature of the method employed here is that if one starts on a hyperbolic periodic orbit (with no nett rotations of its stable and unstable manifolds) then $\eta_s$ will never give $\omega(\eta_s,\xi_s)=0$ with $\lambda(\eta_s)<0$.  In contrast, starting on an elliptic periodic orbit the method should give non-existence fairly fast (look back at the simple example of the pendulum).  Thus it seems $P_0$ is an unfortunate choice of symmetry plane.  From  Figure~\ref{fig:sos}, $\theta=\pi$ looks more hopeful, though even there not all the island chains have elliptic points. Indeed, looking at the symmetry plane $P_\pi$ for $\mu=0.1$ (actually, $P_0$ for $\mu=0.9$, which is equivalent and saved further code changes) in Figure~\ref{fig:lambdapi}, we still see blue regions around some of the resonances.
\item The method does not necessarily eliminate points with unbounded orbits; for $\mu=0$, these lie in $L^2>2r$. This is because it addresses non-existence of any invariant surface through the given orbit transverse to $\xi$, including unbounded surfaces, not just tori.
\end{enumerate}

\begin{figure}[htb]
\centering
\includegraphics[width=4in]{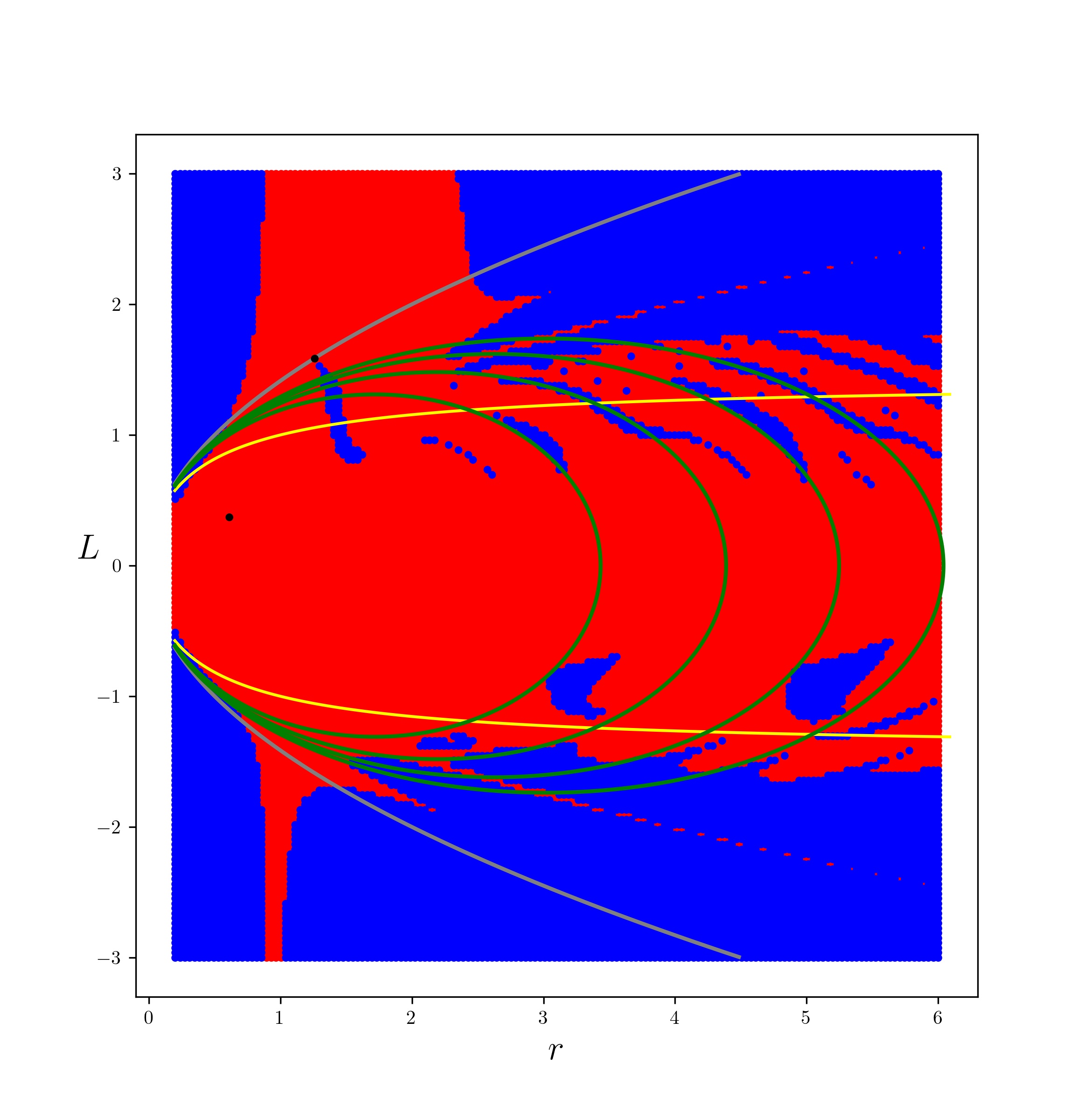}
\caption{Converse KAM results using the symmetric formulation (the rest of the setup same as in Figure~\ref{fig:lambda}).}
\label{fig:eta}
\end{figure}

Next we test the symmetric formulation exploiting the time-reversal symmetry, for which the results are shown in Figure~\ref{fig:eta}. Here the singularities for $\eta$ are just the two Lagrange points near the second body. Comparing with Figure~\ref{fig:lambda}, we see that the symmetric formulation gives more nonexistence results, and so yields a  clearer picture than the general one. This is especially evident around the resonances. There seems to be a general agreement, but the symmetric formulation revealed an extended nonexistence region near $r=1$ for large positive $L$ as well as a much smaller one towards large negative $L$, both lying outside the boundary of $\Sigma$ (black curves in Figure~\ref{fig:lambda}). This difference could be due to the possibility that $\xi$ becomes parallel to $V$. As discussed in sections \ref{subsec:nonexistence}-\ref{subsec:symmetric}, in those places the general formulation will not detect nonexistence, staying consistent with the transversality requirement of the tori class under investigation, whereas the symmetric formulation is not restricted by this and justly will.

\begin{figure}[htb]
\centering
\includegraphics[width=4in]{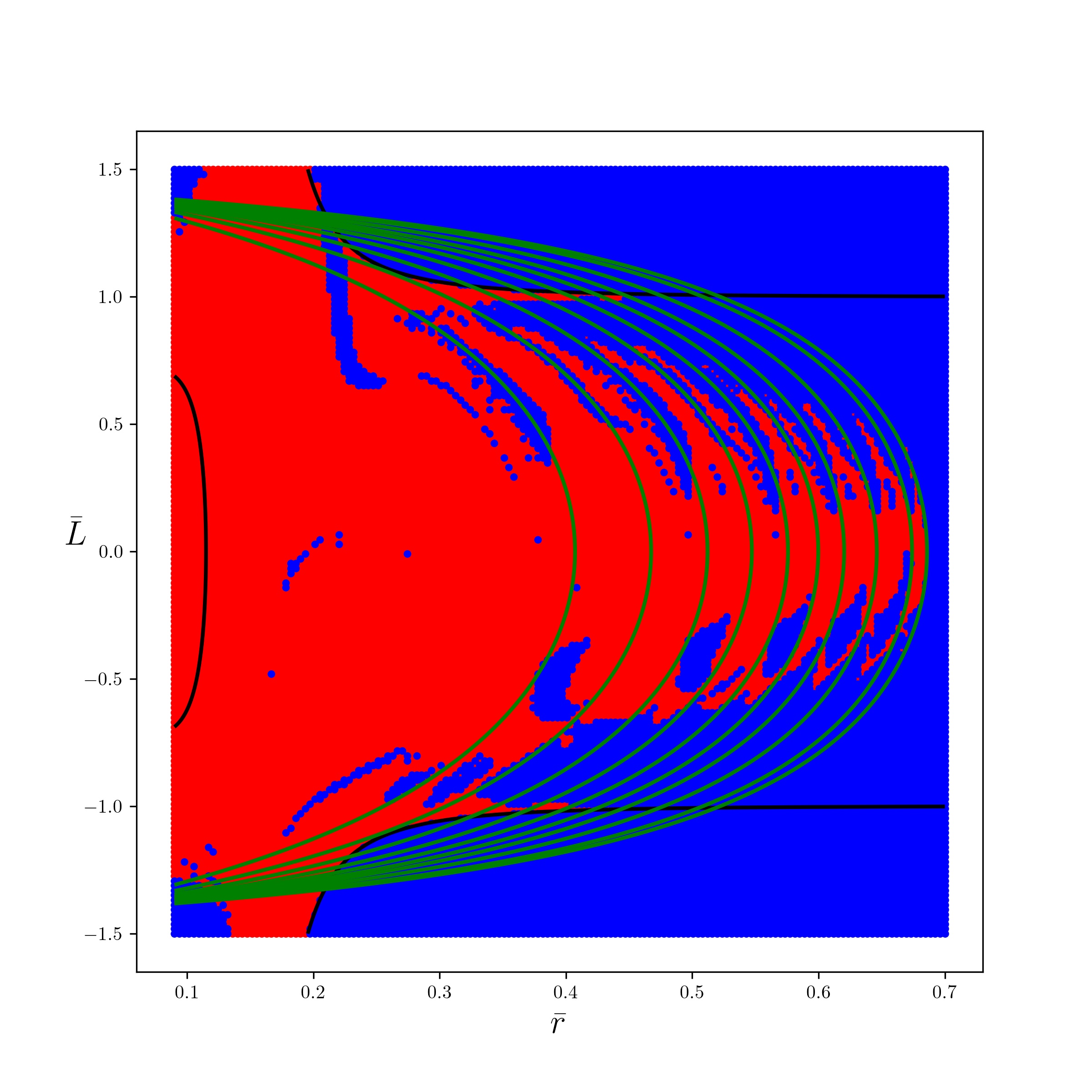}\vspace{-0.25cm}
\caption{Figure~\ref{fig:lambda} replotted in $(\bar{r},\bar{L})$-coordinates, and resonances (green) with winding ratios $n/4$ for $n=9,13,17,21,25,29,33,39,47,51$.}
\label{fig:bars}
\end{figure}

For the PCR3BP, better coordinates to plot the converse KAM results on the symmetry planes might be $\bar{L} = L/\sqrt{r}$ and $\bar{r} = r/(r+m)$ for some $m$, because then the non-escape region for $\mu=0$ is $|\bar{L}|<\sqrt{2}$, the circular orbits are on $|\bar{L}|=1$ and $\bar{r}$ turns $r=(0,\infty)$ into the bounded interval $\bar{r}=(0,1)$.  Figure~\ref{fig:bars} shows the same results as before replotted using coordinates $(\bar{r},\bar{L})$ for $m=5$ instead of $(r,L)$. 
 
\begin{figure}[hp]\vspace{-1cm}
\centering
\includegraphics[width=4.7in]{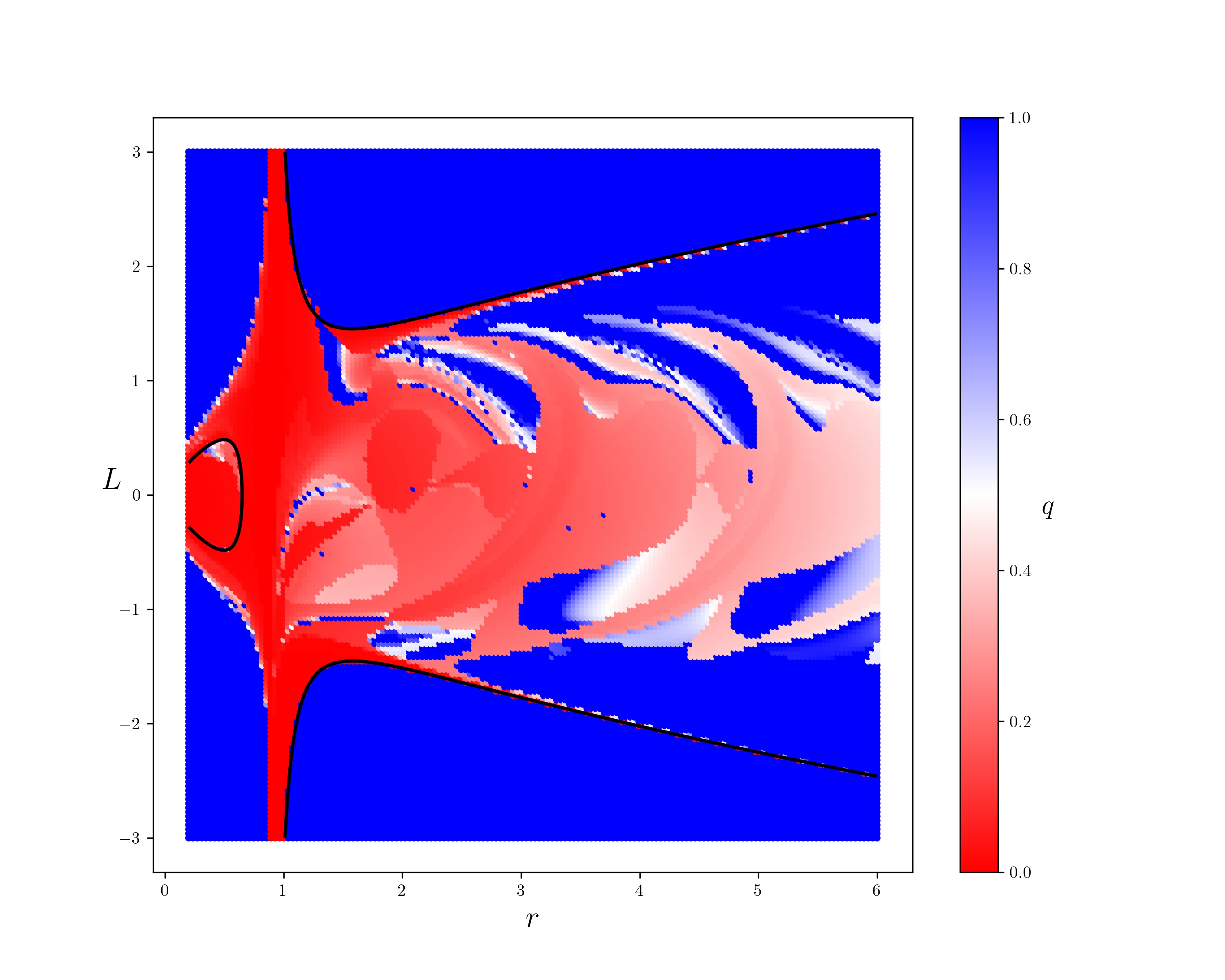}\vspace{-0.6cm}
\caption{Converse KAM measurement. Figure~\ref{fig:lambda} replotted with hues ranging from fast detection of nonexistence (deep red) to slowest one (light blue) and ultimately to no detection at all (deep blue) within timeout.}
\label{fig:speed}
\end{figure}

Finally, Figure~\ref{fig:speed} shows where the method works faster. Here we use $q=1-t_{\text{r}}/t_{\text{out}}$ as a measure, where $t_{\text{r}}$ is the remaining time after the converse KAM condition was satisfied up to timeout $t_{\text{out}}$. As before, blue indicates that nonexistence was not detected until $t_{\text{out}}$, but now the nonexistence region is coloured according to $q$, with a darker red indicating shorter times and light blue longer ones. As we see, nonexistence is detected faster near the two bodies and around the strip $r=1$, and slower for larger distances.

\begin{figure}[hp]\vspace{-0.5cm}
\centering
\includegraphics[width=4.7in]{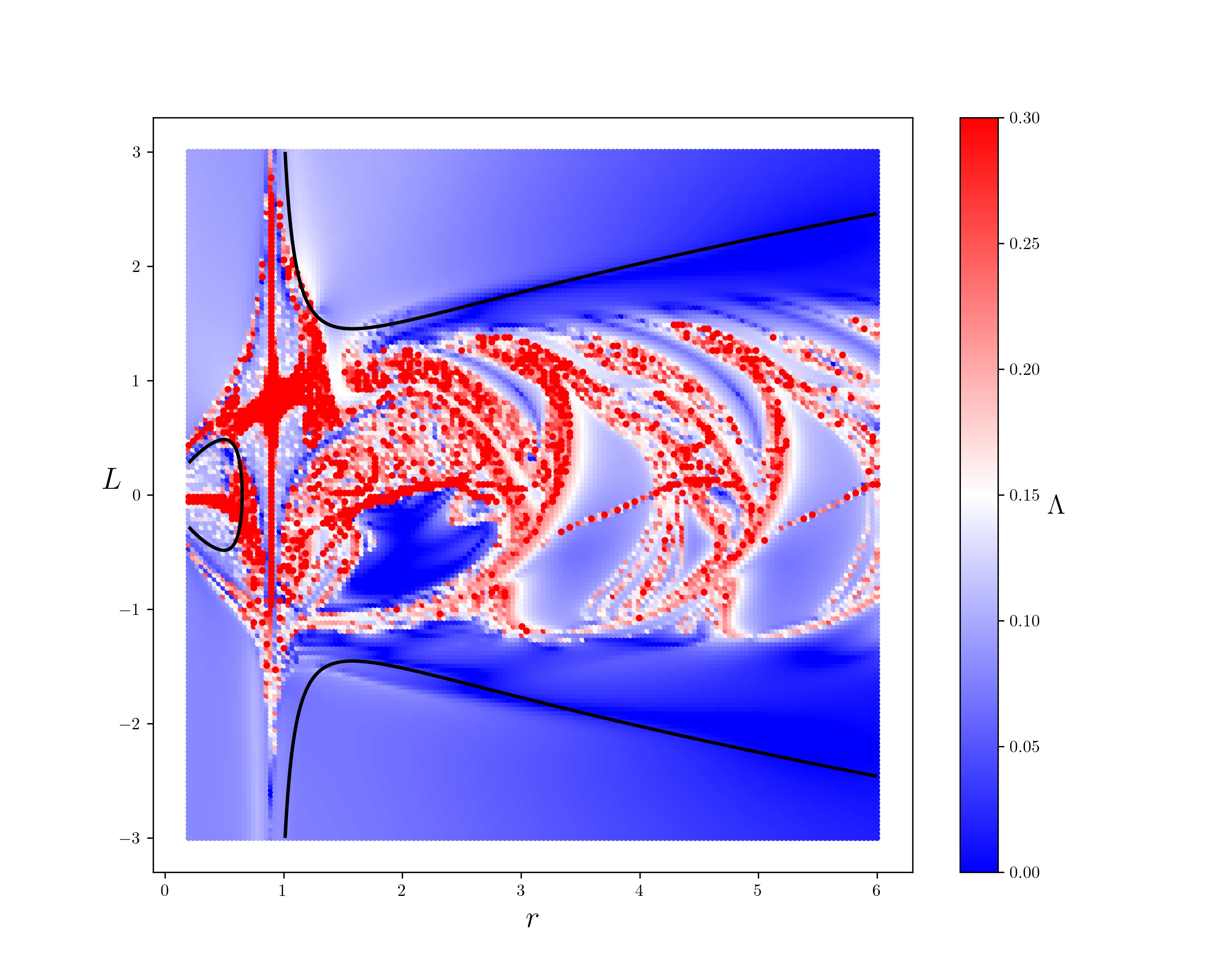}\vspace{-0.6cm}
\caption{The Lyapunov exponent $\Lambda$ for the same setup as in Figure~\ref{fig:speed}.}
\label{fig:lyapunov}
\end{figure}

As a comparison, we also computed the Lyapunov exponent $\Lambda$ as a chaos indicator. Figure~\ref{fig:lyapunov} shows the simple estimate $\Lambda = \frac{1}{t_{\text{out}}}\log\left(|\xi_{t_{\text{out}}}|/|\xi_0|\right)$ for the same initial conditions and parameters as in Figures~\ref{fig:lambda} and \ref{fig:speed}. Initial conditions with red hues appear to have significantly positive Lyapunov exponent. The majority of the orbits though have rather small values (blue hues) and seem to need further investigation (e.g., longer times, check convergence, etc.) to decide if they are chaotic or not. Compared to converse KAM, indeed we see some of them lying inside the nonexistence region (perhaps indicating invariant tori of a different class), but most of them lie in the blue inconclusive region of Figure~\ref{fig:speed}. However, we do see some good agreement near $r=1$ and quite good agreement around the upper tongues near the resonances.

We have recently developed a refinement of Lyapunov exponent calculation in the Hamiltonian context, to distinguish more clearly between positive and zero. In particular, it is expected to distinguish more clearly the zero exponent that arises for trajectories on invariant tori with a smooth conjugacy to a constant vector field.  We will report on that elsewhere.

\section{Improvements to make}
This was an implementation of the method of \cite{M18} on a significant test problem with a non-trivial foliation.  It has demonstrated that the method is usable and useful.  Nevertheless, there are many improvements we could make.

Firstly, other choices of foliations might be better. For example, to treat initial conditions inside the orbit of the secondary (question 2 of the Introduction), it would be better to base the foliation on the primary instead of the centre of mass. Or, to investigate invariant tori around the secondary requires a foliation adapted to it (question 3 of the Introduction). In principle, one could choose a foliation (with singularities) that is simultaneously adapted to all regions, cf.~\cite{DM}. In general, however, choosing a foliation (and a $\lambda$ accordingly) that is suitable for the system under investigation and satisfies the method's requirements can be challenging. A systematic way of constructing the foliation or the direction $\xi$ would be desirable. Based on the ideas of \cite{DM}, the vector field $\xi$ (\ref{eq:xi}) and the corresponding $\lambda$ used here offer one way out, but it would be good to improve and explore other means; this is currently under investigation. 

Secondly, it would be best to regularise collisions with the two heavy bodies (in particular, to compute more accurately without adaptive time-step the trajectories that pass close to the secondary).  This is relatively straightforward (e.g.~the Thiele-Burrau regularisation \cite{Sz}) and avoids the need to adapt timestep for close approaches. To make a first test of the method, however, we decided not to implement this.

Thirdly, we could do with a better surface of section.  Theoretically, one can make a transverse section by finding the continuation of the circular periodic orbits and deforming $p_r=0$ and its boundary $r=L^2$ to make the boundary be the continued circular orbits.  This is Birkhoff's prescription, but it is not explicit.  An alternative is to adopt the procedure of \cite{DW}.  We note also that it seems insufficient to study initial conditions on symmetry planes; unlike in simple problems like the standard map or the H\'enon map, there does not appear to be a dominant symmetry.  To take care of this, we should apply the ``killends'' extension of the method, as outlined in \cite{M18}.

Fourthly, it would probably be better to use a symplectic integrator, to respect the symplectic structure of the problem.  An example is the Stormer-Verlet method, which works for mechanical systems in a rotating frame \cite{M92}. However, the converse KAM method works fairly fast, so that there is no need to go to longer integration times where standard methods like the Runge-Kutta scheme used here might fail in accuracy.

Fifthly, it would be good to devise an escape condition.  This would eliminate many initial conditions with $K\ge 0$ from being on invariant tori of any class, for example.

Last but not least, it would be good to extend the method to higher DoF so that we could treat the planar elliptic or the spatial circular or the spatial elliptic restricted three-body problem, or even the general three-body problem.  A paper on this is in preparation. The idea is to restrict attention to Lagrangian submanifolds transverse to a given Lagrangian foliation.  All tori constructed by usual KAM proofs are Lagrangian, so it is appropriate to restrict to Lagrangian submanifolds. They are also $C^1$ graphs of actions as functions of angles, thus transverse to the foliation by the surfaces of fixed angles, which are Lagrangian. The tangent plane to a Lagrangian submanifold is Lagrangian.  There is a cyclic partial order on Lagrange planes at a given point and the dynamics preserves this order. Thus if there is an invariant Lagrangian submanifold transverse to the Lagrangian foliation, the trajectory of the tangent plane to the foliation can not cross that of the tangent plane to the submanifold. This provides a sufficient condition for non-existence of such a submanifold.

\section{Conclusion}
We have applied a method to establish regions of phase space through which pass no invariant tori of given class, to the planar circular restricted  three-body problem.  It finds large regions of non-existence of tori, which mainly appear to correspond to trajectories whose orbit crosses that of the secondary. We also detect non-existence from some resonances, but not all, which appears to be because in the plane where we chose initial conditions they happen to have hyperbolic periodic points rather than elliptic or inversion hyperbolic.  Nonetheless, the method gives significant restrictions on the regions where stable orbits for a planet could orbit a binary star.  

The study has indicated various issues with applying the method and suggested improvements for the future.

We anticipate the method being useful in many applications.  Specifically, we plan to apply it to magnetic fieldline flow and to guiding-centre motion in magnetic fields.  The method could also find applications to classical models of chemical reaction dynamics.

\section*{Acknowledgements}

This work was partly supported by a grant from the Simons Foundation (601970, RSM).

\newpage
\bibliography{NoToribibfile}

\end{document}